\documentclass[a4paper,11pt]{amsart}
\usepackage[table,dvipsnames]{xcolor}
\usepackage[latin1]{inputenc}
\usepackage{amsmath, amsthm, amssymb, amsfonts, amscd}
\usepackage[english]{babel}
\usepackage[all]{xy}
\usepackage[ruled,vlined]{algorithm2e}
\usepackage[linktocpage=true]{hyperref}
\usepackage{color}
\usepackage[normalem]{ulem}
\usepackage[table]{xcolor}
\usepackage{pinlabel}
\usepackage{tikz-cd}
 \usepackage{booktabs}
\usepackage{enumitem}
\usepackage{graphicx}
\usepackage{subcaption}
\usepackage[margin=3.15cm]{geometry}
\usepackage[font=small,labelfont=bf]{caption}

\graphicspath{ {figures/} }

\title{Hypergraphs for multiscale cycles\\ in structured data}

\author[Barbensi, Yoon, Madsen, Ajayi, Stumpf, Harrington]{Agnese Barbensi, Hee Rhang Yoon, Christian Degnbol Madsen, Deborah O. Ajayi, Michael P.H. Stumpf, Heather A. Harrington}

\begin{document}

\maketitle

\begin{abstract}
Scientific data has been growing in both size and complexity across the modern physical, engineering, life and social sciences. Spatial structure, for example, is a hallmark of many of the most important real-world complex systems, but its analysis is fraught with statistical challenges. Topological data analysis can provide a powerful computational window on complex systems. Here we present a framework to extend and interpret persistent homology summaries to analyse spatial data across multiple scales. We introduce hyperTDA, a topological pipeline that unifies local (\textit{e.g.}~geodesic) and global (\textit{e.g.}~Euclidean) metrics without losing spatial information, even in the presence of noise. Homology generators offer an elegant and flexible description of spatial structures and can capture the information computed by persistent homology in an interpretable way. Here the information computed by persistent homology is transformed into a weighted hypergraph, where hyperedges correspond to homology generators. We consider different choices of generators (\textit{e.g.}~matroid or minimal) and find that centrality and community detection are robust to either choice. We compare hyperTDA to existing geometric measures and validate its robustness to noise. We demonstrate the power of computing higher-order topological structures on spatial curves arising frequently in ecology, biophysics, and biology, but also in high-dimensional financial datasets. We find that hyperTDA can select between synthetic trajectories from the landmark 2020 AnDi challenge and quantifies movements of different animal species, even when data is limited. 
\end{abstract}

\section*{Introduction}
Intricate spatial structures are a hallmark of many complex and natural systems. The movement of individuals through space is one example of a natural system where understanding the spatial (or spatio-temporal) structures is key for distilling mechanistic insights from complex data \cite{nathan2022big}.
Since complex systems tend to produce complicated and highly correlated data, a purely descriptive exploratory analysis offers only limited insights. Model-based analysis, by contrast, suffers from (i) a lack of reliably useful mathematical or computer models; and (ii) the statistical challenges in fitting such models against data and in choosing the most appropriate model~\cite{kirk2013model}. 
\par
The analysis of models and data with topological data analysis (TDA) is gaining increasing traction. Persistent homology (PH), a prominent tool in TDA, computes a geometric summary of the structure and shape of complex systems, including spatio-temporal mathematical models \cite{edelsbrunner2010computational, ghrist2008barcodes, wasserman2018topological}. Persistent homology builds a collection of discrete approximations on the data called {\em simplicial complexes}, computes topological invariants, such as homology, which provides a summary of topological features (\textit{e.g.}~connected components, loops, and voids). 
Crucially, the PH algorithm involves the construction of a so-called filtration of simplicial complexes from the data, which quantifies the persistence of topological features in the data across changes in scale and resolution \cite{edelsbrunner2010computational, ghrist2008barcodes, roadmap}. In this way, a threshold parameter value often required in data analysis is sidestepped.

\par
\begin{figure}[!ht]
    \centering
    \includegraphics[width = 15cm]{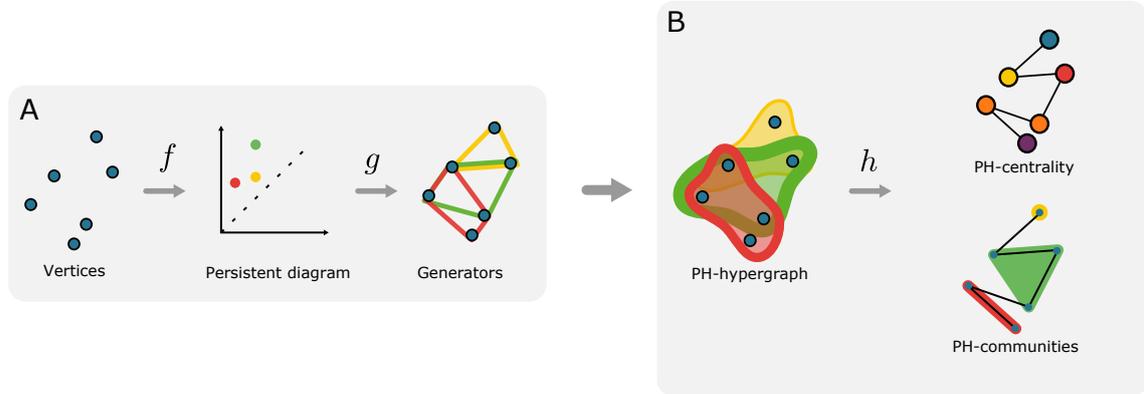}
    \captionsetup{margin=-0.8cm}

    \caption{\textbf{Schematic of hyperTDA pipeline} \textbf{(A)} Given a collection of vertices, we first   
    compute the persistence diagram. The persistence diagram depends on a choice of filtration $f$ to build a nested sequence of simplicial complexes. We choose a way of computing generators $g$ and then compute them for each point in the persistence diagram \textbf{(B)} We encode the information present in the persistence diagram and generators by constructing a hypergraph whose nodes are the vertices and whose hyperedges are the generators. The hyperedges are weighted by the persistence of the corresponding generators. We analyse and interpret the higher-order relationships in the data by performing network theoretic methods, given by h. Specifically, we compute PH-centrality, which quantifies the intensity of interactions between vertices, and PH-community detection, which partitions the vertices into communities that share homology generators.}
    \label{fig:pipeline}
\end{figure}

Persistent homology provides a multiscale analysis of connectivity inherent in data. The powerful algebraic representation that encodes the topological features in persistent homology can be visualised by a \textit{persistence diagram} (see Figure~\ref{fig:pipeline}(A)). 
A persistence diagram is a multiset (a generalisation of the conventional set where elements can appear more than once). Each $(x,y)$ point in the persistence diagram corresponds to a topological feature in the dataset; the importance of a topological feature is measured by its \textit{persistence}, \textit{i.e.}~the norm $|y-x|$ of the corresponding point in the multiset. Much statistical machinery has been developed to vectorise persistence diagrams for statistical or machine learning \cite{pun2022persistent, townsend2020representation}. Recent work has studied persistence diagrams with simplicial complexes \cite{jaramillo2022barcode}. PH analysis has led to parameter inference, pattern detection and non-trivial associations in highly structured and complex data \cite{hiraoka2016hierarchical, mcguirl2020topological, stolz2022multiscale, thorne2022topological,townsend2020representation,  vipond2021multiparameter}.

\par
Since its inception some 20 years ago, topological data analysis and persistent homology have proven to be (often vastly) more informative than conventional summary statistics of complex data. However, it is often more important to identify the likely cause of certain topological features than to merely characterise the full topology. Recent progress in extracting \textit{generators} of persistence diagrams, \textit{i.e.}~vertices in the dataset that form a \textit{cycle} corresponding to a specific $(x,y)$ point in the persistence diagram, has inspired new research that utilises generators as geometric and spatial features of data \cite{benjamin2022homology, emmett2015multiscale, henselmanghristl6, kovacev2016using, LiMinimalCycle}. For example, generators can localise the specific subchain responsible for the entanglement that distinguishes knotted and unknotted homologous proteins \cite{benjamin2022homology}. But different choice of vertices for the same topological feature can lead to different descriptions and interpretations of data. 
\par
Here, we extend topological data analysis from one generator to all generators. To overcome the potential complications arising from the non-uniqueness of generators we develop two complementary approaches: matroid generators, which are efficiently computed using matroid theory \cite{henselmanghristl6}, and minimal generators \cite{emmett2015multiscale, LiMinimalCycle}, which we tailor to sequential datasets. We encode each generator as a hyperedge which forms the \textit{PH-hypergraph}. The PH-hypergraph encodes both local (geodesic) and global (Euclidean) geometric information of the data \cite{emmett2015multiscale,  LiMinimalCycle}.

The PH-hypergraph can be analysed using the armoury of modern graph theory and network science. Rather than applying PH to study network data \cite{aktasPersistenceHomologyNetworks2019a, FengSpatial}, we use network theory to analyse persistent homology. We show that hypergraph centrality measures and community detection distil information and insights from persistent homology generators. We refer to the full pipeline as \textit{hyperTDA}, as it combines hypergraph theory and topological data analysis (TDA). We show that hyperTDA can detect higher-order topological structures that are not captured by conventional TDA; moreover, it is robust to noisy data and to different choices of generators. 
\par
We showcase hyperTDA on a wide range of datasets and highlight its ability to detect and quantify aspects of complex spatio-temporal data.
Motivated by protein carbon backbone \cite{benjamin2022homology} and animal movement trajectories \cite{nathan2022big}, here we analyse curves in three-dimensional space. Each curve can be translated into data as a collection of ordered vertices in 3-space and the number of vertices is the length of the curve. 
We first consider two synthetic datasets: (i) two ``hand-drawn" curves to illustrate the concepts presented 
and then 800 \textit{random curves} of four different lengths. 
Our next dataset, called the \textit{AnDi trajectories}, consists of curves generated by five anomalous diffusion models \cite{munoz2021objective, munoz2020anomalous}. Lastly, we examine \textit{movement tracks}, which is a collection of experimental movement trajectories from nematodes moving in an agar plate and from zebrafish embryos migrating towards a wound \cite{garriga2016expectation, liepe2016accurate}. Tracks of nematodes consist of vertices in two dimensions; therefore, we convert the data into vertices in three dimensions by encoding the time variable in the $z$-coordinate.

\section*{From persistence diagrams to hypergraphs}
Our goal is to systematically analyse the multiscale cycles present in the data. We compute persistent homology, encode all generators by creating a PH-hypergraph, and then analyse it with network theory methods.

\subsection*{From data to generators} We use persistent homology to compute a multiscale summary of cycles present in a curve~\cite{benjamin2022homology, emmett2015multiscale, stolz2022multiscale}. Given a spatial curve (\textit{i.e.}~a sequence $\mathcal{V}$ = $\{v_1 \ldots v_n\}$ of ordered vertices, where $v_i$ are points in 3-space), we choose a filtration $f$ of simplicial complexes that represents the connectivity among vertices at various distance parameters. Here, $f$ is the Vietoris-Rips filtration \cite{ghrist2008barcodes}, which is based on the pairwise Euclidian distances of the vertices $v_i$. We then compute the dimension-1 PH, which summarises the one-dimensional topological features (\textit{i.e.}~cycles) present in the data. The birth and death of cycles across the filtration are summarised in a persistence diagram (Figure~\ref{fig:pipeline}(A)). Each point in the persistence diagram represents a cycle among the vertices, with $x$ and $y$ coordinates corresponding to \textit{birth} and \textit{death} parameters of the filtration. The life span of each point (\textit{birth} parameter $-$ \textit{death} parameter) is called its \textit{persistence} and provides a measure of its significance~\cite{ghrist2008barcodes, zomorodian2005computing}.

Given a point in the persistence diagram, the corresponding cycle is represented by a collection of vertices called a \textit{generator}. Generators are not unique and are therefore often not analysed, even though they encode how different vertices are organised to form cyclic features. Here, we present two methods $g$ for assigning generators to a persistence diagram: the \textit{matroid generators} and the \textit{minimal generators}~(see~\nameref{sec:codeAvail} and Supplementary Information). Matroid generators are the generators returned by the Eirene software \cite{henselmanghristl6} (which relies on matroid theory). While the matroid generators are fast to obtain, the selection of generators is not canonical or optimal. In order to compute generators in a principled way, we designed the minimal generators to capture the geometry and structure of the underlying curve. The minimal generators are computed by first finding the length-minimising generators \cite{LiMinimalCycle} and then minimising the number of jumps between non-consecutive vertices. The matroid and minimal generators are similar to each other (Figure \ref{fig:validation1}(F)), and the hyperTDA outputs are highly correlated (Figure \ref{fig:validation1}(G)). Since the results are comparable and the matroid generators are computationally faster to obtain, we present the results using the matroid generators. 

\subsection*{From generators to PH-hypergraphs}
Once a choice of generators $g$ is made, all generators of a persistence diagram are translated into a \textit{hypergraph}. A hypergraph is a generalisation of a graph, in which hyperedges are arbitrary non-empty subsets of nodes \cite{bick2021higher}. More precisely, a weighted hypergraph is given by a triple $(\mathcal{N},\mathcal{E},\mathcal{W})$, where $\mathcal{E}$ is a collection of non-empty subsets of $\mathcal{N}$, with hyperedge-weights $\mathcal{W} \subset \mathbb{R}$. We construct the  \textit{PH-hypergraph} $\mathcal{H} = (\mathcal{V},\mathcal{G}, \mathcal{P})$, where each hyperedge $g_i \in \mathcal{G} =  \{g_1, \cdots g_h \}$ corresponds to the vertices in a generator. The weight $p_i \in \mathcal{P} = \{p_1, \cdots p_h \}$ is the persistence of the corresponding point in the persistent diagram (Figure~\ref{fig:pipeline}(B)).

To interpret the structural information of the curve's shape provided by the PH-hypergraph, we require a method $h$ to analyse it. We perform two analyses. First, we compute the vertex centrality~\cite{tudisco2021node}, which returns a vector whose entries give the importance of each vertex in the curve. This vector of vertex values -- \textit{PH-centrality} -- provides a ranking of vertices in the hypergraph based on their hyperedge membership and weights (Figure~\ref{fig:pipeline}(B)). In other words, the importance of a vertex depends on the importance of its connections; so a vertex with high centrality is a member of generators with large persistence. The second method we employ to analyse the PH-hypergraph is community detection to find functional modules in the curve. Community detection methods find densely-connected groups of vertices, with sparse connections between groups (relative to what one would expect at random) \cite{python_louvain}. Each vertex $v_i$ is assigned to a unique community $c_j$ where $j = 1, \ldots, k$ and $k$ is the number of communities. The communities are encoded in a matrix of size $n \times k$, where the ${(i,j)}$-entry is $1$ if vertex $v_i$ belongs to community $c_j$ and $0$ otherwise. We refer to the partition matrix of $\mathcal{H}$ as the \textit{PH-community matrix}, where each column is a \textit{PH-community} that corresponds to a functional module induced by higher-order interactions among the vertices in the curve (Figure~\ref{fig:pipeline}(B)).

\section*{Interpretation and robustness of hyperTDA}
HyperTDA detects and quantifies spatial properties of data by encoding higher-order interactions that are often undetectable with standard methods. Here we validate and interpret the output of hyperTDA on the random curves dataset. Curves in this dataset have various lengths and are generated as self-avoiding random walks. Further, we demonstrate that hyperTDA is robust to perturbation of data and different choices of generators.

\subsection*{Validation and correlation of hyperTDA with geometric invariants}
\begin{figure}[!ht]
\centering
\hspace*{-15pt}
    \includegraphics[width = 16cm]{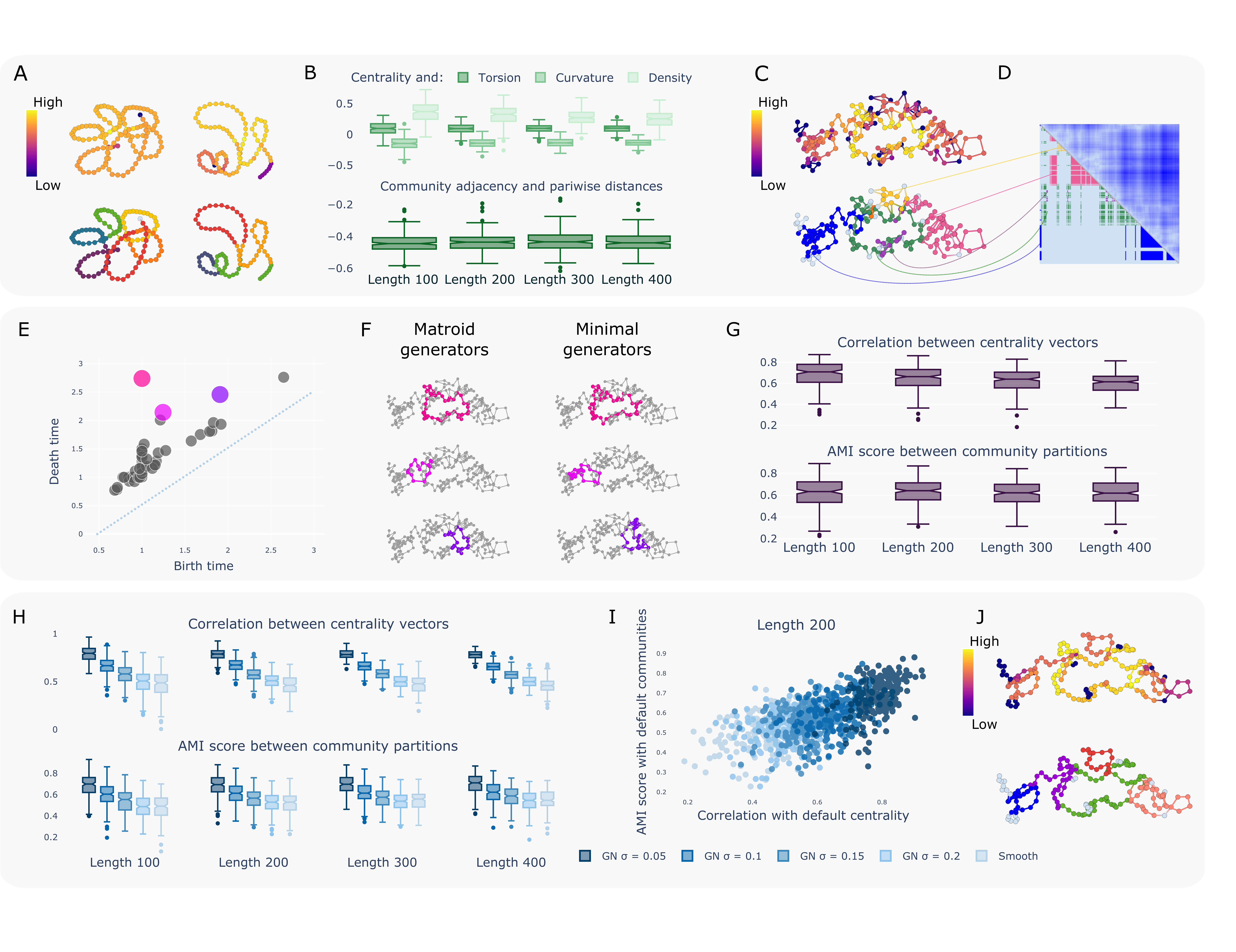}
    \captionsetup{margin=-0.8cm}
    \caption{\textbf{HyperTDA quantifies spatial complexity and is robust to noise.} \textbf{(A)} Two illustrative curves with vertices coloured according to PH-centrality (top) and to PH-communities (bottom) \textbf{(B)} (Top) Correlation between PH-centrality and curvature, torsion, and density vectors for random curves. (Bottom) Correlation between PH-community matrices and pairwise distance matrices for random curves. \textbf{(C)} Example random curve $C$ of length 200. The colour of the vertices represents centrality values (top) and community membership (bottom, light grey vertices are singleton).  \textbf{(D)} The upper triangle shows the pairwise-distance matrix of $C$, with lower values corresponding to lighter shades of blue. The lower triangle shows the PH-community matrix corresponding to $\mathcal{H}(C)$, where non-zero entries are coloured according to community membership. \textbf{(E)} The persistence diagram of $C$, with example points highlighted. \textbf{(F)} The matroid generators and minimal generators for the highlighted points in (E). \textbf{(G)} (Top) Correlation between PH-centrality from matroid generators and from minimal generators. (Bottom) Adjusted mutual information (AMI) score between the partition induced by PH-communities from matroid generators and from minimal generators. \textbf{(H)} Consistency of PH-centrality and PH-communities among random curves and perturbed/smooth curves; top row shows the correlation among PH-centrality and the bottom row shows AMI score between PH-community matrices. \textbf{(I)} AMI score plotted against correlation, for random curves of length 200. \textbf{(J)} The curve obtained by smoothing $C$. Vertices are coloured according to PH-centrality (top) and PH-communities (bottom). PH-centrality identifies the central loop as the predominant feature, and the PH-community matrix identifies the same loop as a community.}
    \label{fig:validation1}
\end{figure}

We first apply the hyperTDA pipeline on two illustrative curves (Figure~\ref{fig:validation1}(A), left, flower and right, roller-coaster). In the flower curve, PH-centrality (top) is essentially homogeneous, thus adequately recognising each of the petals as equally important features; however, the PH-communities (bottom) identifies all the loops as distinct communities. Similarly, in the roller-coaster curve, PH-centrality highlights the most significant geometric and structural features of the curve, while the PH-communities detect the individual loops.

To extensively test and validate the effectiveness of hyperTDA, we apply the pipeline to a dataset of 800 randomly generated piecewise linear open curves in 3-space. The curves are four different lengths, given by 100, 200, 300 and 400 vertices in the curve.  We construct the PH-hypergraph and compute the PH-centrality values and the PH-community matrices. Figure~\ref{fig:validation1}(C) shows an example curve of length 200, with vertices coloured according to PH-centrality values and community membership, respectively. PH-centrality identifies the loop in the middle of the structure as the curve's dominant feature; PH-communities partition the vertices into six distinct modules, one of which (green vertices) corresponds to the central loop (vertices in light blue in Figure~\ref{fig:validation1}(C), bottom row, are singletons that do not belong to any community).

The pairwise distance matrix of a curve has entries equal to the pairwise distances between vertices, which coarsely describes the curve's 3D structure, but fails to encode higher-order features. For example, the dominant loop in Figure~\ref{fig:validation1}(C) detected by hyperTDA cannot be inferred from the pairwise distance matrix, see Figure~\ref{fig:validation1}(D). Recall that the PH-community matrix is the matrix whose ${(i,j)}$-entry is 1 if the \textit{i}$^\text{th}$ and \textit{j}$^\text{th}$ vertices belong to the same community and 0 otherwise; see Figure~\ref{fig:validation1}(D). We compute the correlation between the pairwise distance and PH-community matrices. As shown in Figure~\ref{fig:validation1}(B), we observe a clear negative correlation, with median values of approximately $-$0.5. The negative correlation indicates that the PH-communities recover the global spatial structure coming from the Euclidean metric while providing additional information (\textit{e.g.}~functional modules given by higher-order interactions, as the loop in Figure~\ref{fig:validation1}(C)). Indeed, each community represents clusters of higher-order interactions among vertices (\textit{e.g.}~loops) that are often overlooked by the Euclidean metric alone. 

Similarly, PH-centrality quantifies the data in a way that complements Euclidean distances. We compute the correlation between PH-centrality and torsion, curvature, and density (computed point-wise as the fraction of vertices with distance less than 2). As shown in Figure~\ref{fig:validation1}(B, top), we notice a very weak positive and negative correlation with torsion and curvature (median of approximately 0.1 and -0.14, respectively), and a slightly stronger positive correlation with density (median of approximately 0.27). 

\subsection*{HyperTDA outputs consistent results across different choice of generators}
We find that the minimal generators and the matroid generators are very similar, as shown in Figure~\ref{fig:validation1}(F). Both PH-centrality and PH-communities computed from minimal generators and the matroid generators are highly correlated as shown in Figure~\ref{fig:validation1}(G). To compute correlations of the partitions induced by PH-communities, we used adjusted mutual information (AMI) score.

\subsection*{HyperTDA is robust to noisy data}

One of fundamental advantages of PH is its robustness to noise, which is guaranteed by theoretical results~\cite{cohen2007stability}. These results, however, do not extend automatically to generators. To assess the robustness of generators to noisy data, we apply hyperTDA to the dataset obtained by adding random Gaussian noise of increasing intensity to each vertex (up to $\delta = 0.2$, which is a fifth of the distance between adjacent vertices). We call the resulting curves \textit{perturbed curves}. Further, we create \textit{smooth curves} via a smoothing algorithm. Smoothing strongly affects the point cloud positions while maintaining a reasonable amount of the curve's dominant features; see Figure~\ref{fig:validation1}(J). To check that the curve's overall structure is still recovered, we compute the AMI score between the PH-community matrix of the random and perturbed curves. We also compute the correlations between centralities of the random and perturbed curves. Both PH-centrality and PH-community matrices present high similarities (Figure~\ref{fig:validation1}(H-J)).

\subsection*{PH-communities encode Euclidean and geodesic metrics}
\begin{figure}[!ht]
    \hspace{-7mm}
    \includegraphics[width = 14cm]{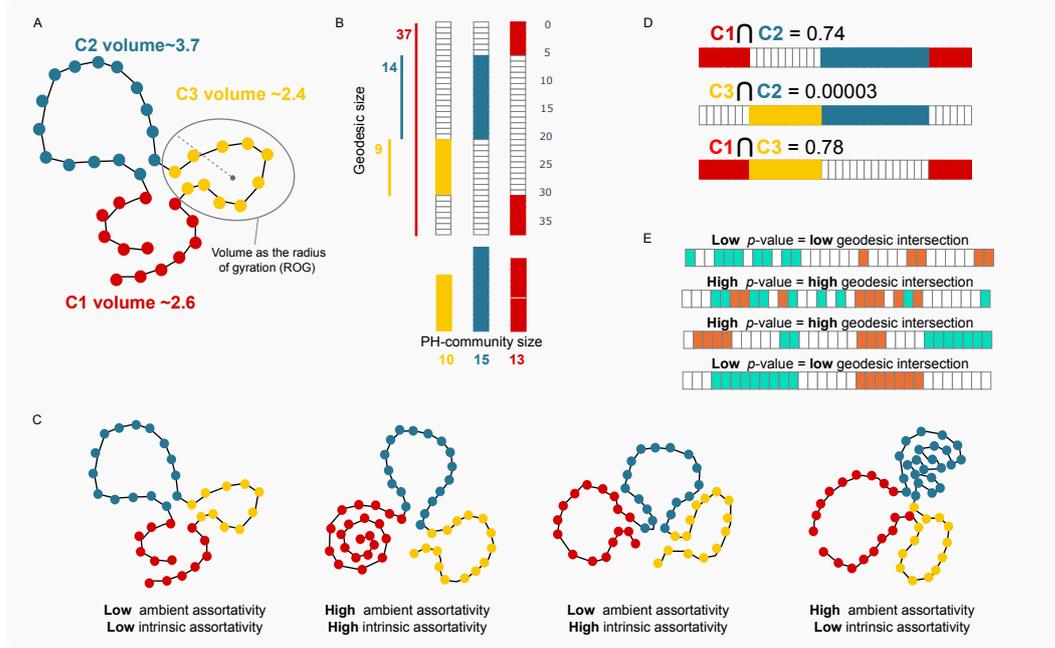}
    \captionsetup{margin=-0.8cm}
    \caption{\textbf{Euclidean and geodesic analysis of PH-communities}  \textbf{(A)} A curve with vertices coloured according to PH-community membership. \textbf{(B)} The panel shows a schematic representation of each community vector, together with their sizes and geodesic sizes. \textbf{(C)} Four examples of curves exhibiting different ambient and intrinsic assortativites \textbf{(D)} Schematic of pairwise geodesic intersection between communities, with examples of pairs with low and high intersection values. \textbf{(E)} Pairwise geodesic intersection values for the communities C1, C2 and C3 in (A).}
    \label{fig:commanalysis}
\end{figure}

We investigate how well PH-communities align with the ambient data space (\textit{i.e.}~Euclidean metric) and the intrinsic metric (\textit{i.e.}~geodesic, which in this case is relative to a discrete space curve). We compute three quantities associated with the communities. First, we estimate the \textit{volume} filled by a community as the radius of gyration of the sub-curve spanned by its vertices. Then, we define the \textit{geodesic size} of a community as $v_{min} - v_{max}$, where $v_{min}$ is the first vertex in the ordered vertex set and $v_{max}$ is the largest vertex in the ordered vertex set. Lastly, we compute the \textit{size} of a PH-community as the number of its vertices (Figure \ref{fig:commanalysis}(A-B)).

In network science, a partition is said to exhibit \textit{assortativity} if similar vertices are more likely to be in the same partition than vertices in different groups. We measure similarity of vertices using either the ambient or intrinsic metric. 

\begin{enumerate}

\item We say that a PH-hypergraph exhibits \textit{ambient assortativity} if vertices in the same PH-community are likely to be very close in space. We define the ambient assortativity as $ \frac{1}{k} \cdot \sum_{i=1}^k  \frac{\text{Size}(c_i)}{ \text{Volume}(c_i)} $, the average ratios between sizes and volumes of all non-singletons communities $c_1, \dots c_k$ (Figure~\ref{fig:commanalysis}(C)). Curves with highly packed communities will then have high ambient assortativity (Figure \ref{fig:commanalysis}(C)). 

\item Similarly, we say that a PH-hypergraph exhibits \textit{intrinsic assortativity} if vertices in the same PH-community are likely to be separated by a small geodesic size. We measure the intrinsic assortativity by $\frac{1}{k} \cdot \sum_{i=1}^k  \frac{\text{Size}(c_i)}{\text{Geodesic size}(c_i)}$, the average ratios between sizes and geodesic sizes of all communities $c_1, \dots c_k$ (Figure~\ref{fig:commanalysis}(C)). With this definition, curves with communities of small geodesic sizes have high intrinsic assortativity. (Figure \ref{fig:commanalysis}(C)). 

\end{enumerate}

By construction, PH-communities define clusters of vertices highly interacting with each other. From a geodesic point of view; however, pairs of communities might appear more or less overlapping, depending on the curve's structure. For example, in Figure~\ref{fig:commanalysis}(D), C2 and C3 are completely separated. In contrast, C1 and C3, and C1 and C2 intersect geodesically. To quantify this phenomenon, we define the \textit{geodesic intersection} as the \textit{p}-value of the Mann-Whitney U-test between the sets of vertices in the communities. Choosing the \textit{p}-value is needed to avoid artefacts due to noise and isolated vertices. Figure~\ref{fig:commanalysis}(E) shows examples of different pairs of PH-communities with high and low intersection values.

\section*{Results}
We first demonstrate the ability of hyperTDA to distinguish diffusion behaviours. Specifically, we show that PH-communities and PH-centrality present statistically significant differences among different diffusion models. We then showcase hyperTDA in the analysis of species movements from experimental data.

\subsection*{HyperTDA distinguishes and selects anomalous diffusion models}
Anomalous diffusion is the transport behaviour of a particle in which the mean square displacement~(MSD) is proportional to~$t^{\alpha}$ for time $t$ and some parameter $\alpha \neq 1$. Anomalous diffusion phenomena have been observed across several fields, including molecular dynamics~\cite{di2018anomalous, krapf2015mechanisms,sabri2020elucidating}, ecology~\cite{humphries2012foraging}, and finance~\cite{plerou2000economic}. The ubiquitous presence of anomalous diffusion has inspired important research efforts to investigate and model its mechanisms; these efforts resulted in several different anomalous diffusion models~(ADMs).

\begin{figure}[!ht]
    \hspace{-7mm}
    \includegraphics[width = 14cm]{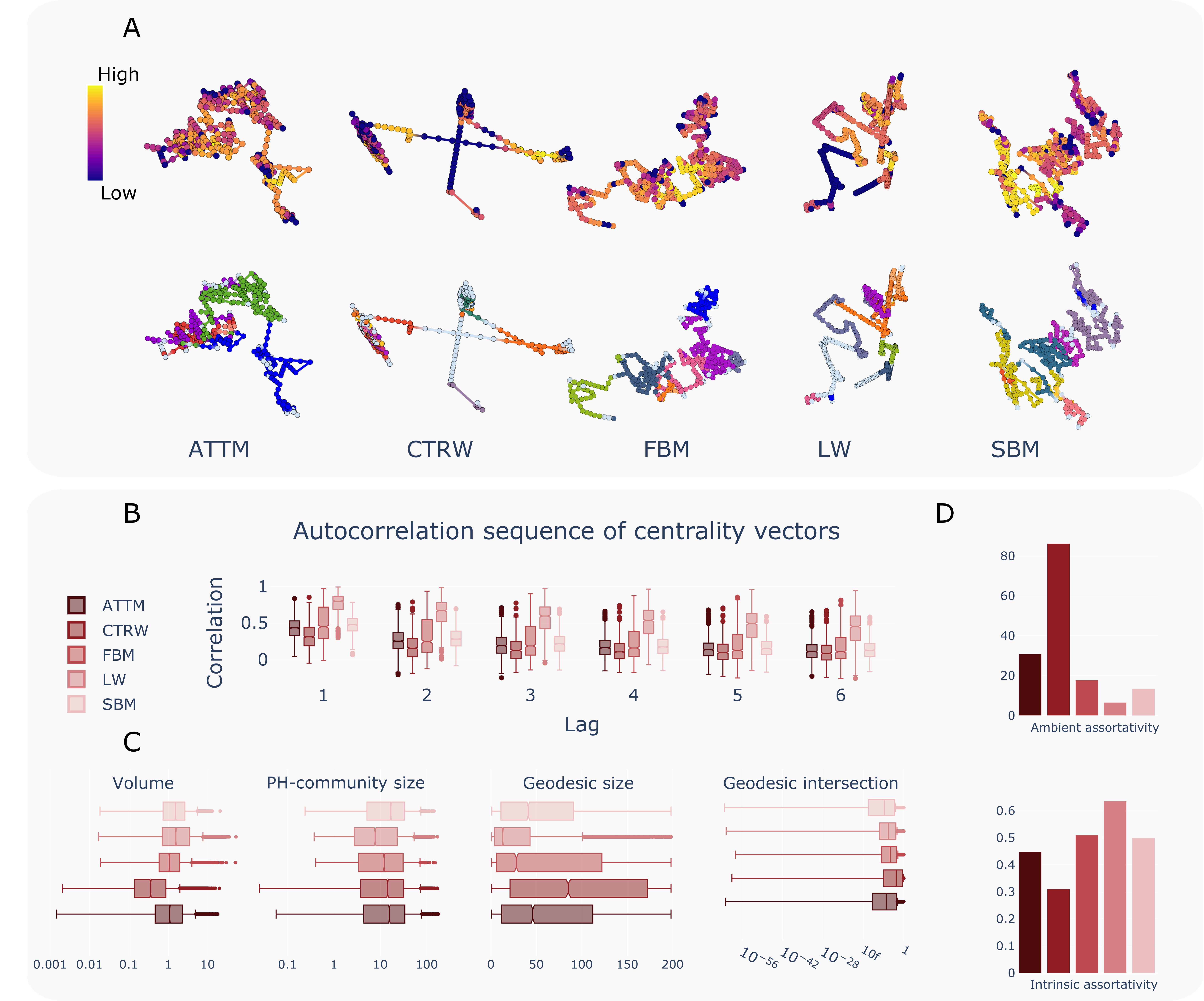}
    \captionsetup{margin=-0.8cm}
    \caption{\textbf{HyperTDA distinguishes anomalous diffusion models (ADMs)}  \textbf{(A)} A representative trajectory for each model, in the following order: ATTM, CTRW, FBM, LW, SBM.  On the top row, vertices are coloured by PH-centrality. On the bottom, vertices are coloured by community membership. \textbf{(B)} Auto-correlation of centrality distributions, by model; for each value $m$ in $\{1,\cdots{}, 6\}$ the boxplots show the distribution of correlation values between centrality vectors and their values shifted by a lag $m$. \textbf{(C)} Analysis of the PH-communities, by model. Box plots show the distributions of volume, size, geodesic size and geodesic intersection are shown on the right. \textbf{(D)} Bar plots showing by model ambient and intrinsic assortativities. }
    \label{fig:diffusioncommunity}
\end{figure}

Inferring the correct model for experimental trajectories is a difficult problem~\cite{munoz2021objective}. One reason making the model detection challenging is that the 3D structure of a trajectory is strongly influenced by its MSD (and thus it depends on the parameter $\alpha$). On the other hand, each ADM is defined for a range of values of $\alpha$ and thus produces trajectories exhibiting a range of different MSDs. We apply the hyperTDA pipeline to the anomalous random walks analysed in the 2020 AnDi challenge~\cite{munoz2021objective, munoz2020anomalous}, which has been the first synergistic effort to develop methods able to classify individual trajectories among several models of diffusion. These models are the annealed transient time motion~(ATTM), the continuous time random walk~(CTRW), the fractional Brownian motion~(FBM), the L\'evy walk~(LW), and the scaled Brownian motion~(SBM) (Figure~\ref{fig:diffusioncommunity}(A)). These ADMs have been selected by the challenge organisers due to their biological relevance~\cite{munoz2021objective}, as they model diffusion processes in living cells (CTRW, ATTM and FBM \cite{ krapf2015mechanisms, saxton2007biological, verdier2021learning}), animal foraging (LW \cite{humphries2012foraging, verdier2021learning}), and human white matter (SBM \cite{fieremans2016vivo}).

We apply hyperTDA to trajectories arising from these five models. We aim to detect and interpret model-specific differences arising from hyperTDA analysis of trajectories. We simulate 1000 length-200 trajectories for each of the five ADMs and apply the hyperTDA pipeline (compute PH-centrality and PH-communities) to all trajectories. We find the pattern of PH-centrality values depends on the underlying model and detects model-specific differences (Figure~\ref{fig:diffusioncommunity}(B)). 
Next we compare the different models by computing the PH-community size, volume, geodesic size, and pairwise geodesic intersections 
(Figure~\ref{fig:diffusioncommunity}(C)). Independent of the scaling exponent~$\alpha$, the distribution of these PH-community analysis gives a statistically significant difference between models. 

We use network theory to further quantify the model-specific structural features highlighted by the PH-communities. We define the \textit{model ambient assortativity} as the average of the ratios between PH-community sizes and volumes over all non-singleton communities in the model. Similarly, the \textit{model intrinsic assortativity} is the average of the ratios between PH-community sizes and geodesic sizes over all non-singleton communities in the model (Figure~\ref{fig:diffusioncommunity}(D)). Assortativity provides a refined interpretation of the model-specific structure. For example, CTRW trajectories exhibit the highest discrepancy between Euclidean and geodesic metrics, as they score the lowest intrinsic assortativity (Figure~\ref{fig:diffusioncommunity}(D)) and highest ambient assortativity. Thus, their structure is characterised by compact, highly interconnected communities that are spread out along the curve. We highlight that standard computations of PH are blind to these features (see Supplementary Information, Figure \ref{fig:landscapes}).

\begin{figure}[!ht]
    \centering
    \includegraphics[width = 13cm]{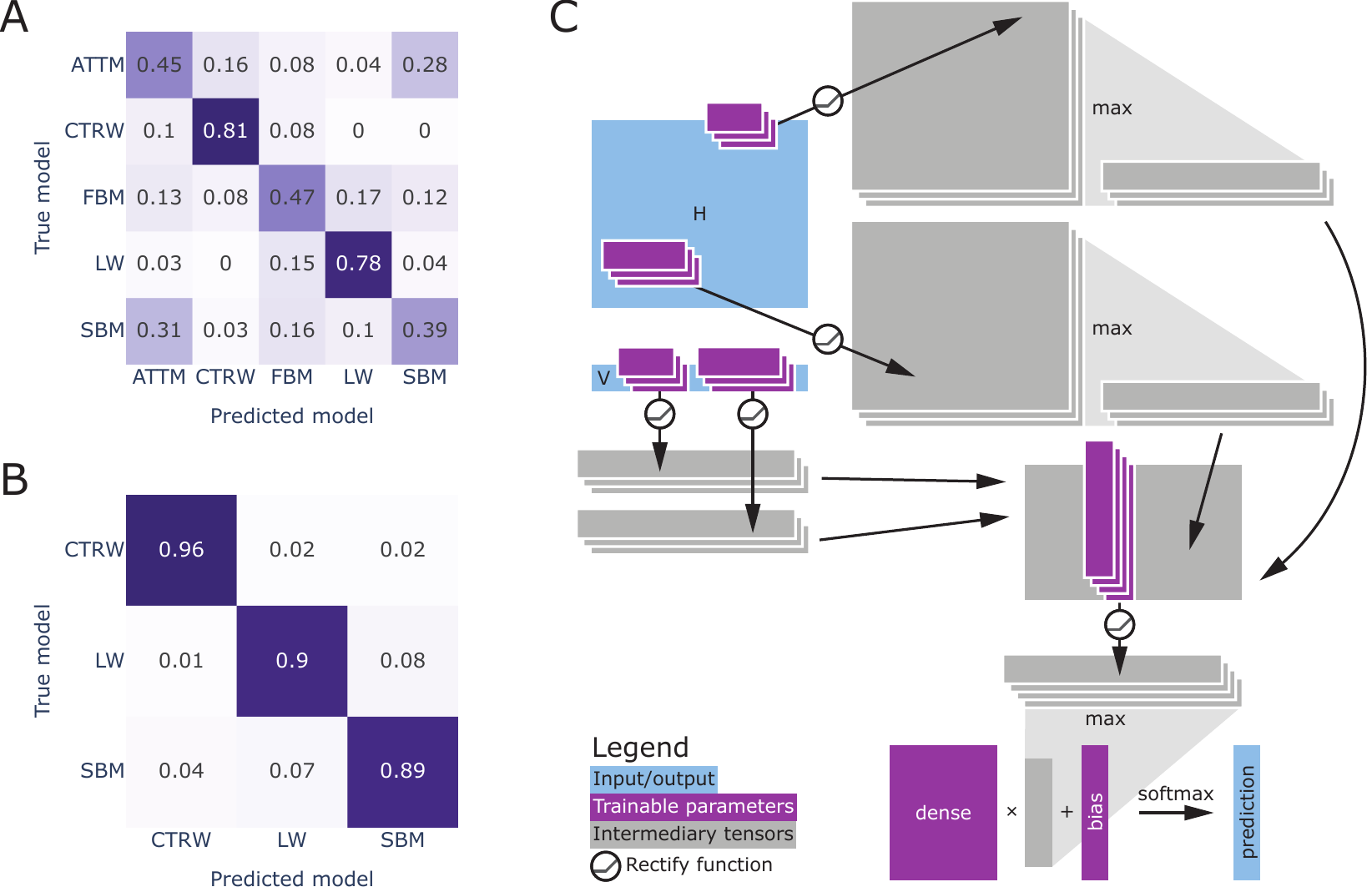}
    \captionsetup{margin=-0.8cm}

    \caption{\textbf{HyperTDA + CNN successfully select ADMs.} \textbf{(A)} The confusion matrix summarising the classification accuracy between all the 5 models. \textbf{(B)} The confusion matrix summarising the classification accuracy between CTRW, LW and SBM models. \textbf{(C)} A schematic description of the CNN architecture given a single datapoint input, that is: a matrix H with dimensions $\#\mathrm{hyperedges}\times\#\mathrm{vertices}$, indicating hyperedge memberships, and vertex centrality values V. Arrows follow the flow of the forward pass. Rectangular boxes represent matrices and vectors, and stacking indicates a third dimension. All trainable parameters are convolutional except when otherwise indicated. The label ``max'' stands for \textit{maxpooling}. For simplicity, the schematic only shows two convolutional kernels with 3 filters in the first layers, and 4 filters in the last.
    }
    \label{fig:CNN}
\end{figure}

Given trajectory data, we next explore whether hyperTDA can predict the model that generated the data. To test this hypothesis, we construct a convolutional neural network (CNN) that takes as input PH-hyperedges and node centrality values. A schematic of the architecture is shown in Figure~\ref{fig:CNN}(C), and the classification results are summarised in Figure~\ref{fig:CNN}(AB). Despite some difficulties in distinguishing the ATM trajectories from FBM and SBM, we obtain an accuracy comparable to the performance of ranked participants of the AnDi challenge~\cite{munoz2021objective, munoz2020anomalous}. We emphasise that the classification relies only on the structural description of the trajectories, effectively disregarding any temporal information. Further, when we restrict to a comparison among CTRW, LW, and SBM, we achieve high accuracy of approximately 92\%.

\subsection*{HyperTDA distinguishes species movements}

The increased availability of technologies to record animal trajectories has inspired research efforts aimed at inferring information on animal behaviour from their tracked movements~\cite{garriga2016expectation,lewis2021learning,maekawa2020deep, nathan2022big}. These approaches have the potential to reveal insight into interactions between individuals, and between animals and their environment, with unprecedented details~\cite{nathan2022big}. One of the main challenges in movement research is to find efficient and unsupervised ways to fragment animal trajectories into units known as \textit{behavioural nodes}~\cite{nathan2008movement}. Similarly, the recent advances in cell imaging have been stimulating progress in the study of cell migration~\cite{liepe2016accurate,masuzzo2016taking,svensson2018untangling,yamada2019mechanisms}. Natural challenges in this context include quantifying homogeneity of  migration patterns in a population and the identification of biases and persistence (\textit{i.e.}~the time period during which a cell moves in the same direction) in the underlying random walks~\cite{beltman2009analysing,krause2014steering,svensson2018untangling}. 

\subsubsection*{PH-centrality identifies behavioural nodes in nematode trajectories}
\begin{figure}[!ht]
    \includegraphics[width = 14cm]{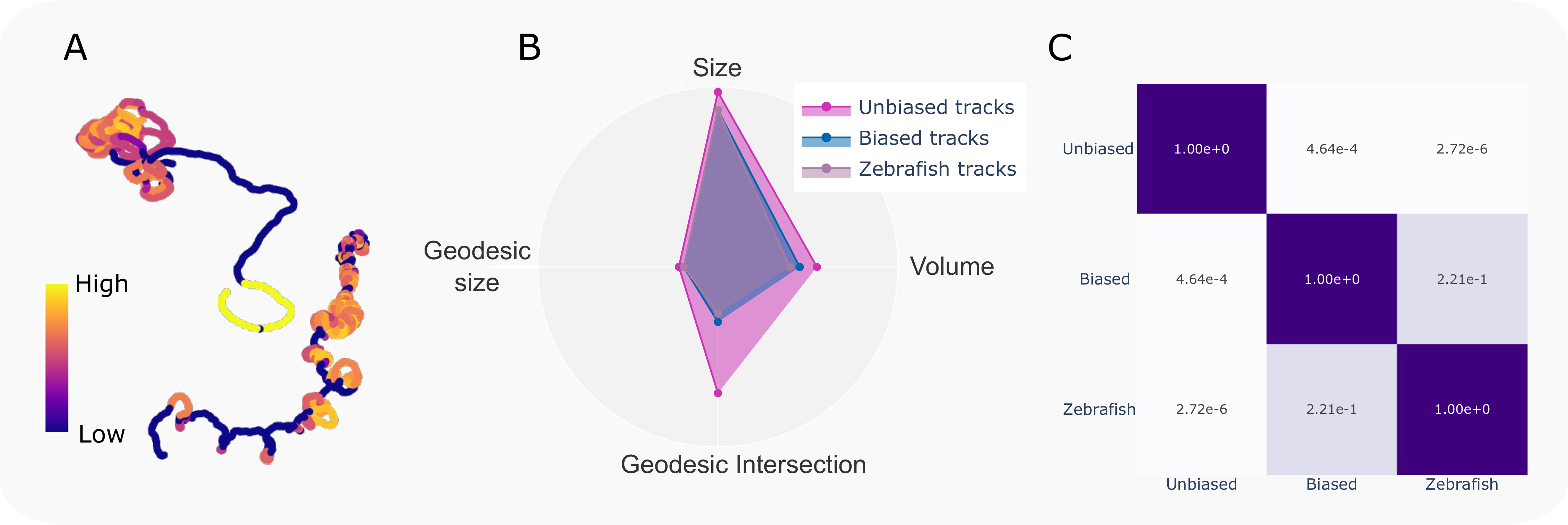}
    \captionsetup{margin=-0.8cm}

    \caption{\textbf{Animal movement and cell migration trajectories} \textbf{(A)} Two 2-D trajectories, with vertices plotted by centrality value. Centrality distinguishes different behaviours in terms of intensity of local search, looping behaviour and relocation. \textbf{(B)} Radar chart given by the medians of volume, size, geodesic size and pairwise geodesic intersection distributions for the unbiased and biased simulated data, and for the zebrafish trajectories. Results are shown with minimal generators. \textbf{(C)} Heat-map of pairwise Kolmogorov-Smirnov test on the cumulative distributions obtained by concatenating geodesic intersection, volume, size and geodesic size distributions. Results are shown with minimal generators.}
    \label{fig:animalcell}
\end{figure}

We apply hyperTDA to animal movement trajectory data. The first dataset consists of six distinct solitary nematodes moving in an agar plate \cite{garriga2016expectation}. PH-centrality clearly distinguishes different behaviours of the nematode (Figure \ref{fig:animalcell}(A)). Indeed, portions of the curve assigned low PH-centrality (Figure \ref{fig:animalcell}(A), blue) can be interpreted as a behaviour known as relocation. Conversely, regions of the curve with increased PH-centrality correspond to behaviours known as local search and looping; furthermore, higher PH-centrality values correspond to more localised and intense search and looping behaviour (Figure \ref{fig:animalcell}(A), yellow). These behavioural classification of vertices obtained from PH-centrality is consistent with those of previous studies~\cite{garriga2016expectation}, yet it is unsupervised and arguably more nuanced. We find that PH-communities identify individual loops (and possibly groups of intertwined loops) of the worm's trajectory as functional modules (see Supplementary Information). 

\subsubsection*{HyperTDA reveals migration bias of zebrafish embryos in \textit{in vivo} tracks}
As a final application, we apply the pipeline to assess the bias in neutrophil cell tracks extracted from zebrafish embryos migrating towards a wound~\cite{liepe2016accurate}. Neutrophils, a type of white blood cell, tend to migrate toward wounds; thus, the expectation is that they follow a \textit{biased persistent random walk}~\cite{liepe2016accurate}. To quantify the presence of the bias, we perform a statistical comparison of PH-communities in zebrafish tracks with those of biased and unbiased persistent random walks on an ellipse; all data is taken from~\cite{liepe2016accurate}. All simulated trajectories have length 20 while the lengths of \textit{in vivo} tracks vary. We thus partition the latter into segments of length 20 before applying hyperTDA; results are shown in Figure~\ref{fig:animalcell}(BC). PH-communities of biased and unbiased simulated trajectories present statistically significant differences (Figure~\ref{fig:animalcell}(BC)). Indeed, as expected, the PH-communities of biased tracks are on average more compact, less overlapping, and have smaller geodesic size than the unbiased tracks. Further, we observe a striking similarity between the \textit{in vivo} and the biased trajectories (Figure~\ref{fig:animalcell}(BC)). Based on the pairwise Kolmogorov-Smirnov tests, we reject the hypothesis that trajectories in zebrafish are unbiased. Lastly, the ambient and intrinsic assortativities are higher for zebrafish and biased trajectories than the unbiased trajectories (see Supplementary Information). All combined, the analysis demonstrates the bias induced by approaching the wound in the experimental trajectories.

Cell tracks often lie on curved surfaces~\cite{liepe2016accurate}; for this reason, classical methods necessitate the pre-processing of data via complex manifold learning algorithms. In contrast, the data is not pre-processed and hyperTDA can be applied directly to the dataset. 
 Note that the trajectories analysed are quite short and have very simple structures. Importantly, even in this case, both matroid (see Supplementary Information) and minimal generators (Figure~\ref{fig:animalcell}(BC)) are successful in recognising the differences and common patterns. Minimal generators, in particular, performed better in recognising similarities in biased and \textit{in vivo} trajectories (see also SI, Figures~\ref{fig:animalcellSI} and \ref{fig:animalcellSI2}). This suggests minimal generators as a preferable option when dealing with limited data.

\section*{Discussion}
We propose hyperTDA, a systematic and extensive analysis of persistent homology generators, to interpret complex spatial systems. We demonstrate the effectiveness of hyperTDA in capturing defining features of curves that are not detected with traditional methods. The versatility of  hyperTDA has natural direct applications in different contexts, ranging from the statistical inference of behaviour from species tracks to a nuanced and improved characterisation of polymer structures. We successfully employed this method to distinguish and predict anomalous diffusion models and to identify biases in species tracks. As a further example, we are currently working on an extension of this approach to the study of the 3D organisation of the human genome. 

The pipeline immediately generalises to the study of any point cloud. Naturally, the results and their interpretation will depend on the choice of filtration used to compute PH. While the Vietoris-Rips filtration is a standard choice in persistent homology, different options might be preferable depending on the data. Similarly, variations in methods to analyse the PH-hypergraph allow for a range of possible flavours and interpretations. 

The hyperTDA pipeline involves a choice for computing generators for points of the persistence diagram. We propose minimal generators and matroid generators as two viable options, and we show that the hyperTDA outcomes are strongly correlated. Data is often limited as a consequence of financial constraints or ethical concerns. Even under these conditions, hyperTDA with minimal generators captures the subtle signals revealing common patterns and structural differences. On the other hand, when analysing large and complex datasets, computational efficiency is an issue. In such cases, hyperTDA with matroid generators has the advantage of a fast computation, yet it provides similar hyperTDA results to that of minimal generators. One could adapt hyperTDA to other choices of generators~\cite{ChenGenerators, Dey2019, emmett2015multiscale, EscolarHiraoka, LiMinimalCycle, Obayashi2018VolumeOC, Ripser_involuted} depending on the data and the end goal; clearly the fledgling area of research into, and construction of, topology generators has opened up new mathematical avenues for mathematical exploration; but already they are proving to be immensely useful in topological data analysis. 
\par
HyperTDA is flexible, general, and remarkably and demonstrably robust to noise in data. We can get hopelessly lost in the exploration of high-dimensional data. Here we have demonstrated that hyperTDA allows us to navigate such data, detect and quantify higher-order structures and make sense of complex systems.

\section*{Methods}
\label{sec:Methods}
\subsection*{Computing persistent homology}
Persistent diagrams in dimension one are computed using the Julia software Eirene \cite{henselmanghristl6}. Persistent homology is computed with $\mathbb{F}_2$ coefficients. 

\subsection*{From PH to homology generators}
The first choice of generators is the matroid generators; these are returned by Eirene's default computations  \cite{henselmanghristl6}, which rely on matroid theory. The second choice of generators is called minimal generators which are computed as follows. For each point in the persistence diagram, we compute the length-minimal generator \cite{LiMinimalCycle} with rational coefficients. Nearly all solutions have coefficients in $\{0, \pm 1\}$, so we interpret the solution as having coefficients in $\mathbb{F}_2$. We then perform a jump-minimisation step where we iteratively search for homologous generators with shorter jumps \cite{emmett2015multiscale}. That is, given a generator whose vertices are given by $[ \dots, v_i, v_{i+1} \dots]$, if such generator is homologous to $[ \dots v_i, v_*, v_{i+1} \dots]$ for some $v_i < v_* < v_{i+1}$, then we choose the latter generator. Details can be found in Supplementary Information.

\subsection*{From generators to hypergraphs.}
Given a choice $g$ of generators, we construct the PH-hypergraph as follows. The set of nodes is given by the vertices of the curve. Each homology generator identifies a subset of nodes; for each such subset, we add a hyperedge connecting the corresponding vertices. Then, hyperedges are assigned a weight given by the persistence of the corresponding point in the persistence diagram. 

\subsection*{Analysis of hypergraphs generators}
Centrality is computed using the software from \cite{tudisco2021node}, where a class of spectral centrality measures is defined and analysed. The model in \cite{tudisco2021node} depends on the choice of nonlinear functions, with different choices inducing different centrality \textit{flavours} that dictate how hyperedge weights influence node values. For our analysis, we focus on the \textit{max} centrality flavour, which assigns large values to nodes that are part of at least one \textit{important} edge (see~\nameref{sec:codeAvail}).

As community detection on hypergraphs is computationally challenging, we first flatten each PH-hypergraph into a graph, as follows. The set of nodes is maintained; for each hyperedge involving the subset of nodes N$_i$ we add the complete-subgraph on N$_i$. Once the graph is constructed, we perform community detection using the method implemented in Python's Louvain module \cite{python_louvain} with default parameters (see~\nameref{sec:codeAvail}).

Note that some trajectories in the AnDi and Zebrafish datasets have trivial PH. We removed such trajectories from the analysis. 

\subsection*{Noisy Data} 
To assess the robustness of hyperTDA, we apply the pipeline to successive perturbations of the randomly generated curves of lengths 100, 200, 300 and 400. These perturbations are obtained by adding Gaussian noise to each point in a curve. For each curve we apply the perturbation process 4 times, with intensity $\sigma$ increasing from 0.05 to 0.2 (a fifth of the distance between consecutive vertices in a curve). Further, we obtain smooth versions of each curve by applying the pyknotid \texttt{smooth} function \cite{pyknotid}~(see~\nameref{sec:codeAvail}).

\subsection*{Interpolation and uninterpolation}
The AnDi trajectories have a high variance in distance between consecutive vertices. We thus modify the curves by inserting vertices between consecutive vertices of the original curve before computing hyperTDA. Specifically, we interpolate the length 199 trajectories to make consecutive vertices approximately equidistant. The interpolation is achieved by adding a total of 301 linearly interpolated vertices between adjacent vertices. We first assign a single interpolated vertex to the segment with the largest length. We repeat the process for 301 interpolated vertices. At the end of the process, all the vertices assigned to a given segment are distributed evenly between the endpoints. The resulting interpolated trajectories have lengths of 500. We refer to them as the \textit{interpolated curves}.

The outputs of hyperTDA are reported on the original curve. To do so, before computing PH-community size, geodesic size, and assortativity values, we perform \textit{un}interpolation to obtain centrality values and community membership values on the original curve. This is done by multiplying the PH-community matrix $M$ (of size $500 \times k$, where $k$ is the number of communities) by a sparse matrix $D$ that maps interpolated vertices to \textit{un}interpolated vertices. We choose to construct $D$ with a rolling average weighting, where an interpolated vertex partly maps to both adjacent \textit{un}interpolated vertices, with more weight towards the closest of the two. For more details, see~\nameref{sec:codeAvail}.

\subsection*{HyperTDA analysis}
A detailed explanation of the methods and notebooks containing the relevant code can be found in the hyperTDA GitHub repository~(see~\nameref{sec:codeAvail}).

\section*{}
\subsection*{Acknowledgements}
The authors thank Gregory Henselman-Petrusek and Benjamin Schweinhart for helpful conversations. 

\subsection*{Funding}
AB gratefully acknowledges funding through a MACSYS Centre Development initiative from the School of Mathematics \& Statistics, the Faculty of Science and the Deputy Vice Chancellor Research, University of Melbourne. HY gratefully acknowledges funding through the Mark Foundation for Cancer Research. CDM gratefully acknowledges funding through a University of Melbourne PhD studentship. MPHS is funded through the University of Melbourne DRM initiative, and acknowledges financial support from the Volkswagen Foundation through a ``Life?'' program grant. 
This research was partially supported by The University of Melbourne's Research Computing Services and the Petascale Campus Initiative. This research was undertaken using the LIEF HPC-GPGPU Facility hosted at the University of Melbourne. This Facility was established with the assistance of LIEF Grant LE170100200. HAH gratefully acknowledges funding from EPSRC EP/R005125/1, the Royal Society RGF$\backslash$EA$\backslash$201074 and UF150238. AB, DOA and HAH are grateful to the EPSRC GCRF grant EP/T001968/1, part of the Abram Gannibal Project, for supporting collaborative meetings. HAH, AB and CDM are grateful to the support provided by the UK Centre for Topological Data Analysis EPSRC grant EP/R018472/1. For the  purpose of Open Access, the authors have applied a CC BY public copyright licence to any Author Accepted Manuscript (AAM) version arising from this submission.

\subsection*{Code availability}
\label{sec:codeAvail}
The code and datasets are available online at the project's GitHub repository: \href{https://github.com/degnbol/hyperTDA}{https://github.com/degnbol/hyperTDA}. A complete description of each dataset can be found in the Supplementary Information. The software used to compute minimal generators can be found in the following GitHub repository: \href{https://github.com/irishryoon/minimal_generators_curves}{https://github.com/irishryoon/minimal\_generators\_curves}, and a description of the method is in the Supplementary Information. 

\subsection*{Author contribution statements}
AB, HY, and HAH conceived the research; AB, HY, HAH, MPHS designed the research; AB, HY, DOA, and CDM performed research; AB, HY, and CDM analysed data; and all the authors wrote the paper.


\begin{thebibliography}{10}

\bibitem{aktasPersistenceHomologyNetworks2019a}
Mehmet~E. Aktas, Esra Akbas, and Ahmed~El Fatmaoui.
\newblock Persistence homology of networks: methods and applications.
\newblock {\em Applied Network Science}, 4(1):61, August 2019.

\bibitem{python_louvain}
Thomas Aynaud.
\newblock python-louvain x.y: Louvain algorithm for community detection.
\newblock
  \href{https://github.com/taynaud/python-louvain}{\texttt{https://github.com/taynaud/python-louvain}},
  2020.

\bibitem{beltman2009analysing}
Joost~B Beltman, Athanasius~FM Mar{\'e}e, and Rob~J De~Boer.
\newblock Analysing immune cell migration.
\newblock {\em Nature Reviews Immunology}, 9(11):789--798, 2009.

\bibitem{benjamin2022homology}
Katherine Benjamin, Lamisah Mukta, Gabriel Moryoussef, Christopher Uren,
  Heather~A Harrington, Ulrike Tillmann, and Agnese Barbensi.
\newblock Homology of homologous knotted proteins.
\newblock {\em arXiv preprint arXiv:2201.07709}, 2022.

\bibitem{bick2021higher}
Christian Bick, Elizabeth Gross, Heather~A Harrington, and Michael~T Schaub.
\newblock What are higher-order networks?
\newblock {\em arXiv preprint arXiv:2104.11329}, 2021.

\bibitem{ChenGenerators}
Chao Chen and Daniel Freedman.
\newblock Measuring and computing natural generators for homology groups.
\newblock {\em Computational Geometry}, 43:169--181, 02 2010.

\bibitem{cohen2007stability}
David Cohen-Steiner, Herbert Edelsbrunner, and John Harer.
\newblock Stability of persistence diagrams.
\newblock {\em Discrete \& computational geometry}, 37(1):103--120, 2007.

\bibitem{VinesVineyardsPersistence}
David Cohen-Steiner, Herbert Edelsbrunner, and Dmitriy Morozov.
\newblock Vines and vineyards by updating persistence in linear time.
\newblock In {\em Proceedings of the Twenty-Second Annual Symposium on
  Computational Geometry}, pages 119--126, New York, NY, USA, 2006. Association
  for Computing Machinery.

\bibitem{dabrowski2021topoly}
Pawel Dabrowski-Tumanski, Pawel Rubach, Wanda Niemyska, Bartosz~Ambrozy Gren,
  and Joanna~Ida Sulkowska.
\newblock Topoly: Python package to analyze topology of polymers.
\newblock {\em Briefings in Bioinformatics}, 22(3):bbaa196, 2021.

\bibitem{Dey2019}
Tamal~K Dey, Tao Hou, and Sayan Mandal.
\newblock Persistent 1-cycles: Definition, computation, and its application.
\newblock In {\em International Workshop on Computational Topology in Image
  Context}, pages 123--136. Springer, 2019.

\bibitem{di2018anomalous}
Michele Di~Pierro, Davit~A Potoyan, Peter~G Wolynes, and Jos{\'e}~N Onuchic.
\newblock Anomalous diffusion, spatial coherence, and viscoelasticity from the
  energy landscape of human chromosomes.
\newblock {\em Proceedings of the National Academy of Sciences},
  115(30):7753--7758, 2018.

\bibitem{edelsbrunner2010computational}
Herbert Edelsbrunner and John Harer.
\newblock {\em Computational topology: an introduction}.
\newblock American Mathematical Soc., 2010.

\bibitem{emmett2015multiscale}
Kevin Emmett, Benjamin Schweinhart, and Raul Rabadan.
\newblock Multiscale topology of chromatin folding.
\newblock {\em arXiv preprint arXiv:1511.01426}, 2015.

\bibitem{EscolarHiraoka}
E.G. Escolar and Y.~Hiraoka.
\newblock Optimal cycles for persistent homology via linear programming.
\newblock {\em Optimization in the Real World. Mathematics for Industry}, 13,
  2016.

\bibitem{FengSpatial}
Michelle Feng and Mason~A. Porter.
\newblock Spatial applications of topological data analysis: Cities,
  snowflakes, random structures, and spiders spinning under the influence.
\newblock {\em Phys. Rev. Research}, 2:033426, Sep 2020.

\bibitem{fieremans2016vivo}
Els Fieremans, Lauren~M Burcaw, Hong-Hsi Lee, Gregory Lemberskiy, Jelle
  Veraart, and Dmitry~S Novikov.
\newblock In vivo observation and biophysical interpretation of time-dependent
  diffusion in human white matter.
\newblock {\em Neuroimage}, 129:414--427, 2016.

\bibitem{garriga2016expectation}
Joan Garriga, John~RB Palmer, Aitana Oltra, and Frederic Bartumeus.
\newblock Expectation-maximization binary clustering for behavioural
  annotation.
\newblock {\em PLoS One}, 11(3):e0151984, 2016.

\bibitem{ghrist2008barcodes}
Robert Ghrist.
\newblock Barcodes: the persistent topology of data.
\newblock {\em Bulletin of the American Mathematical Society}, 45(1):61--75,
  2008.

\bibitem{Hatcher}
Allen Hatcher.
\newblock {\em {Algebraic topology}}.
\newblock Cambridge Univ. Press, Cambridge, 2000.

\bibitem{henselmanghristl6}
G.~{Henselman} and R.~{Ghrist}.
\newblock {Matroid Filtrations and Computational Persistent Homology}.
\newblock {\em ArXiv e-prints}, June 2016.

\bibitem{hiraoka2016hierarchical}
Yasuaki Hiraoka, Takenobu Nakamura, Akihiko Hirata, Emerson~G Escolar, Kaname
  Matsue, and Yasumasa Nishiura.
\newblock Hierarchical structures of amorphous solids characterized by
  persistent homology.
\newblock {\em Proceedings of the National Academy of Sciences},
  113(26):7035--7040, 2016.

\bibitem{humphries2012foraging}
Nicolas~E Humphries, Henri Weimerskirch, Nuno Queiroz, Emily~J Southall, and
  David~W Sims.
\newblock Foraging success of biological l{\'e}vy flights recorded in situ.
\newblock {\em Proceedings of the National Academy of Sciences},
  109(19):7169--7174, 2012.

\bibitem{jaramillo2022barcode}
Edgar Jaramillo-Rodriguez.
\newblock Barcode posets: Combinatorial properties and connections.
\newblock {\em arXiv preprint arXiv:2206.05613}, 2022.

\bibitem{kirk2013model}
Paul Kirk, Thomas Thorne, and Michael~PH Stumpf.
\newblock Model selection in systems and synthetic biology.
\newblock {\em Current opinion in biotechnology}, 24(4):767--774, 2013.

\bibitem{kovacev2016using}
Violeta Kovacev-Nikolic, Peter Bubenik, Dragan Nikoli\'{c}, and Giseon Heo.
\newblock Using persistent homology and dynamical distances to analyze protein
  binding.
\newblock {\em Statistical Applications in Genetics and Molecular Biology},
  15(1):19--38, 2016.

\bibitem{krapf2015mechanisms}
Diego Krapf.
\newblock Mechanisms underlying anomalous diffusion in the plasma membrane.
\newblock {\em Current topics in membranes}, 75:167--207, 2015.

\bibitem{krause2014steering}
Matthias Krause and Alexis Gautreau.
\newblock Steering cell migration: lamellipodium dynamics and the regulation of
  directional persistence.
\newblock {\em Nature reviews Molecular cell biology}, 15(9):577--590, 2014.

\bibitem{lewis2021learning}
Mark~A Lewis, William~F Fagan, Marie Auger-Methe, Jacqueline Frair, John~M
  Fryxell, Claudius Gros, Eliezer Gurarie, Susan~D Healy, and Jerod~A Merkle.
\newblock Learning and animal movement.
\newblock {\em Frontiers in Ecology and Evolution}, 9:681704, 2021.

\bibitem{LiMinimalCycle}
Lu~Li, Connor Thompson, Gregory Henselman-Petrusek, Chad Giusti, and Lori
  Ziegelmeier.
\newblock Minimal cycle representatives in persistent homology using linear
  programming: An empirical study with user's guide.
\newblock {\em Frontiers in Artificial Intelligence}, 4, 2021.

\bibitem{liepe2016accurate}
Juliane Liepe, Aaron Sim, Helen Weavers, Laura Ward, Paul Martin, and
  Michael~PH Stumpf.
\newblock Accurate reconstruction of cell and particle tracks from 3d live
  imaging data.
\newblock {\em Cell systems}, 3(1):102--107, 2016.

\bibitem{maekawa2020deep}
Takuya Maekawa, Kazuya Ohara, Yizhe Zhang, Matasaburo Fukutomi, Sakiko
  Matsumoto, Kentarou Matsumura, Hisashi Shidara, Shuhei~J Yamazaki, Ryusuke
  Fujisawa, Kaoru Ide, et~al.
\newblock Deep learning-assisted comparative analysis of animal trajectories
  with deephl.
\newblock {\em Nature communications}, 11(1):1--15, 2020.

\bibitem{masuzzo2016taking}
Paola Masuzzo, Marleen Van~Troys, Christophe Ampe, and Lennart Martens.
\newblock Taking aim at moving targets in computational cell migration.
\newblock {\em Trends in cell biology}, 26(2):88--110, 2016.

\bibitem{mcguirl2020topological}
Melissa~R McGuirl, Alexandria Volkening, and Bj{\"o}rn Sandstede.
\newblock Topological data analysis of zebrafish patterns.
\newblock {\em Proceedings of the National Academy of Sciences},
  117(10):5113--5124, 2020.

\bibitem{munoz2021objective}
Gorka Mu{\~n}oz-Gil, Giovanni Volpe, Miguel~Angel Garcia-March, Erez Aghion,
  Aykut Argun, Chang~Beom Hong, Tom Bland, Stefano Bo, J~Alberto Conejero,
  Nicol{\'a}s Firbas, et~al.
\newblock Objective comparison of methods to decode anomalous diffusion.
\newblock {\em Nature communications}, 12(1):1--16, 2021.

\bibitem{munoz2020anomalous}
Gorka Mu{\~n}oz-Gil, Giovanni Volpe, Miguel~Angel Garcia-March, Ralf Metzler,
  Maciej Lewenstein, and Carlo Manzo.
\newblock The anomalous diffusion challenge: single trajectory characterisation
  as a competition.
\newblock In {\em Emerging Topics in Artificial Intelligence 2020}, volume
  11469, pages 42--51. SPIE, 2020.

\bibitem{nathan2008movement}
Ran Nathan, Wayne~M Getz, Eloy Revilla, Marcel Holyoak, Ronen Kadmon, David
  Saltz, and Peter~E Smouse.
\newblock A movement ecology paradigm for unifying organismal movement
  research.
\newblock {\em Proceedings of the National Academy of Sciences},
  105(49):19052--19059, 2008.

\bibitem{nathan2022big}
Ran Nathan, Christopher~T Monk, Robert Arlinghaus, Timo Adam, Josep Al{\'o}s,
  Michael Assaf, Henrik Baktoft, Christine~E Beardsworth, Michael~G Bertram,
  Allert~I Bijleveld, et~al.
\newblock Big-data approaches lead to an increased understanding of the ecology
  of animal movement.
\newblock {\em Science}, 375(6582):eabg1780, 2022.

\bibitem{Obayashi2018VolumeOC}
Ippei Obayashi.
\newblock Volume optimal cycle: Tightest representative cycle of a generator on
  persistent homology.
\newblock {\em SIAM J. Appl. Algebra Geom.}, 2:508--534, 2018.

\bibitem{roadmap}
Nina Otter, Mason~A Porter, Ulrike Tillmann, Peter Grindrod, and Heather~A
  Harrington.
\newblock A roadmap for the computation of persistent homology.
\newblock {\em EPJ Data Science}, 6:1--38, 2017.

\bibitem{plerou2000economic}
Vasiliki Plerou, Parameswaran Gopikrishnan, Lu{\'\i}s A~Nunes Amaral, Xavier
  Gabaix, and H~Eugene Stanley.
\newblock Economic fluctuations and anomalous diffusion.
\newblock {\em Physical Review E}, 62(3):R3023, 2000.

\bibitem{pun2022persistent}
Chi~Seng Pun, Si~Xian Lee, and Kelin Xia.
\newblock Persistent-homology-based machine learning: a survey and a
  comparative study.
\newblock {\em Artificial Intelligence Review}, pages 1--45, 2022.

\bibitem{sabri2020elucidating}
Adal Sabri, Xinran Xu, Diego Krapf, and Matthias Weiss.
\newblock Elucidating the origin of heterogeneous anomalous diffusion in the
  cytoplasm of mammalian cells.
\newblock {\em Physical Review Letters}, 125(5):058101, 2020.

\bibitem{saxton2007biological}
Michael~J Saxton.
\newblock A biological interpretation of transient anomalous subdiffusion. i.
  qualitative model.
\newblock {\em Biophysical journal}, 92(4):1178--1191, 2007.

\bibitem{scharein2002interactive}
Robert~G Scharein and Kellogg~S Booth.
\newblock Interactive knot theory with knotplot.
\newblock In {\em Multimedia Tools for Communicating Mathematics}, pages
  277--290. Springer, 2002.

\bibitem{stolz2022multiscale}
Bernadette~J Stolz, Jakob Kaeppler, Bostjan Markelc, Franziska Braun, Florian
  Lipsmeier, Ruth~J Muschel, Helen~M Byrne, and Heather~A Harrington.
\newblock Multiscale topology characterizes dynamic tumor vascular networks.
\newblock {\em Science Advances}, 8(23):eabm2456, 2022.

\bibitem{svensson2018untangling}
Carl-Magnus Svensson, Anna Medyukhina, Ivan Belyaev, Naim Al-Zaben, and
  Marc~Thilo Figge.
\newblock Untangling cell tracks: Quantifying cell migration by time lapse
  image data analysis.
\newblock {\em Cytometry Part A}, 93(3):357--370, 2018.

\bibitem{pyknotid}
Alexander~J Taylor and other SPOCK~contributors.
\newblock pyknotid knot identification toolkit.
\newblock \url{https://github.com/SPOCKnots/pyknotid}, 2017.
\newblock Accessed YYYY-MM-DD.

\bibitem{thorne2022topological}
Thomas Thorne, Paul~DW Kirk, and Heather~A Harrington.
\newblock Topological approximate bayesian computation for parameter inference
  of an angiogenesis model.
\newblock {\em Bioinformatics}, 38(9):2529--2535, 2022.

\bibitem{townsend2020representation}
Jacob Townsend, Cassie~Putman Micucci, John~H Hymel, Vasileios Maroulas, and
  Konstantinos~D Vogiatzis.
\newblock Representation of molecular structures with persistent homology for
  machine learning applications in chemistry.
\newblock {\em Nature communications}, 11(1):1--9, 2020.

\bibitem{tudisco2021node}
Francesco Tudisco and Desmond~J Higham.
\newblock Node and edge nonlinear eigenvector centrality for hypergraphs.
\newblock {\em Communications Physics}, 4(1):1--10, 2021.

\bibitem{Ripser_involuted}
Matija \v{C}ufar and \v{Z}iga Virk.
\newblock Fast computation of persistent homology representatives with
  involuted persistent homology, 2021.

\bibitem{verdier2021learning}
Hippolyte Verdier, Maxime Duval, Fran{\c{c}}ois Laurent, Alhassan Cass{\'e},
  Christian~L Vestergaard, and Jean-Baptiste Masson.
\newblock Learning physical properties of anomalous random walks using graph
  neural networks.
\newblock {\em Journal of Physics A: Mathematical and Theoretical},
  54(23):234001, 2021.

\bibitem{vipond2021multiparameter}
Oliver Vipond, Joshua~A Bull, Philip~S Macklin, Ulrike Tillmann, Christopher~W
  Pugh, Helen~M Byrne, and Heather~A Harrington.
\newblock Multiparameter persistent homology landscapes identify immune cell
  spatial patterns in tumors.
\newblock {\em Proceedings of the National Academy of Sciences}, 118(41), 2021.

\bibitem{wasserman2018topological}
Larry Wasserman.
\newblock Topological data analysis.
\newblock {\em Annual Review of Statistics and Its Application}, 5:501--532,
  2018.

\bibitem{yamada2019mechanisms}
Kenneth~M Yamada and Michael Sixt.
\newblock Mechanisms of 3d cell migration.
\newblock {\em Nature Reviews molecular cell biology}, 20(12):738--752, 2019.

\bibitem{zomorodian2005computing}
Afra Zomorodian and Gunnar Carlsson.
\newblock Computing persistent homology.
\newblock {\em Discrete \& Computational Geometry}, 33(2):249--274, 2005.

\end{thebibliography}

\newpage

\appendix
\section{Datasets}

All the dataset we analyse are available in the GitHub repository  and summarised in Table \ref{table:tab1}. 

\subsection{Artificial drawings}
The first dataset of artificial examples consists of two curves generated using the software KnotPlot~\cite{scharein2002interactive}. These curves are available in the \texttt{examples/default\_pipeline} folder in the GitHub repository~(see~\nameref{sec:codeAvail}), together with the software used to analyse them, and a notebook to visualise the results.

\subsection{Random curves}
The second dataset consists of 4 sets of 200 random curves of lengths 100, 200, 300 and 400, respectively. These curves are generated as equilateral self-avoiding random walks (SARW) using the python package Topoly \cite{dabrowski2021topoly}. The distance between two consecutive vertices in a curve is equal to 1. These curves are available in the \texttt{results/Random\_curves} folder in the GitHub repository~(see~\nameref{sec:codeAvail}), together with a notebook to visualise the results. Each curve in this dataset has been perturbed in 5 different ways, listed below. \\

\begin{enumerate}
    \item Adding to each point random Gaussian noise with intensity $\sigma = 0.05$;
    \item Adding to each point random Gaussian noise with intensity $\sigma = 0.1$;
    \item Adding to each point random Gaussian noise with intensity $\sigma = 0.15$;
    \item Adding to each point random Gaussian noise with intensity $\sigma = 0.2$;
    \item Taking a smooth version of the curve by using the \texttt{spacecurves.smooth} function of the \texttt{pyknotid} \cite{pyknotid} Python software.
\end{enumerate}

A notebook showing the perturbation methods can be found in the \texttt{results/Random\_curves} folder in the GitHub repository~(see~\nameref{sec:codeAvail}).

\begin{table}[h!]
\small

\begin{tabular}{||p{2.6cm} |p{2.6cm}| p{2.6cm}| p{2.6cm} | p{2.6cm}||} 
 \hline
\normalsize \textbf{Dataset} & \normalsize  \textbf{Size} & \normalsize  \textbf{Lengths} & \normalsize  \textbf{Interpolated} & \normalsize  \textbf{Perturbed}  \\  [0.9ex] 
\hline\hline
 KnotPlot & 2 & 132, 131 & No & No \\
 \hline
 SARWs & 200 & 100 & No & Yes    \\
 \hline
 SARWs & 200 & 200 & No  & Yes  \\
 \hline
 SARWs & 200 & 300 & No   & Yes \\
 \hline
 SARWs & 200 & 400 & No   & Yes \\
 \hline
 AnDi & 5000 & 199 & Yes to len 500   & No  \\
 \hline
 Nematodes & 6 & 1437-1830  & No  & No   \\
 \hline
 Zebrafish & 163 & 20 & No  & No   \\
 \hline
 BPRWs & 401 & 20 & No   & No  \\ 
 \hline
 PRWs & 401 & 20 & No & No \\
 \hline
\end{tabular}
\caption{A summary of the datasets analysed. }
\label{table:tab1}
\end{table}

\subsection{AnDi}
The third dataset consists of 5000 ADM trajectories, generated using the AnDi challenge's software \cite{munoz2021objective, munoz2020anomalous}. All trajectories have length 199. Each trajectory is generated by first randomly choosing a model between ATTM, CTRW, FBM, LW, SBM, then by randomly choosing the exponents $\alpha$. Exponents are $\leq 2$ and strictly $>0.05$, as smaller exponents produce practically immobile trajectories. Moreover, for CTRW and ATTM $\alpha \leq 1 $, for LW $\alpha \geq 1 $, for FBM $\alpha < 2 $, while for SBM we do not have further restrictions. Each trajectory is interpolated before applying the pipeline. The length of interpolated trajectories is set to 500. For details on the interpolation methods see the Methods section of the main manuscript. The software used to interpolate is available in the GitHub repository. AnDi trajectories are available in the \textit{results/AnDi} folder in the GitHub repository, together with a notebook to visualise the results. Details on the hyperTDA analysis can be also found in the \texttt{examples/with\_interpolation} folder in the GitHub repository~(see~\nameref{sec:codeAvail}). Out of the 5000 trajectories analysed, 53 had trivial PH. 

\subsection{Nematodes}
The fourth dataset consists of the recordings of six distinct solitary nematodes moving in an agar plate, taken from~\cite{garriga2016expectation}. The trajectories are 2-dimensional, and their lengths vary between 1438 and 1830. Before applying the pipeline we added time as a \textit{z}-coordinate. Nematode tracks are available in the \texttt{results/species\_tracks} folder in the GitHub repository~(see~\nameref{sec:codeAvail}), together with a notebook to visualise the results and a notebook that shows the pre-processing. The folder and the notebook include the community partition of each trajectory and the corresponding plots. 

\subsection{Zebrafish tracks}
The fifth dataset consists of 200 length 20 biased persistent random walks (BPRWs) on ellipses, 200 length 20 unbiased persistent random walks (PRWs) on ellipses (taken from \cite{liepe2016accurate}), and 162 length 20 sections of zebrafish embryos \textit{in vivo} tracks. Zebrafish tracks are obtained by cutting length 20 segments from the trajectories available from \cite{liepe2016accurate}. As the average distance between successive points is different between simulated and \textit{in vivo} data, we uniformly scale all the trajectories to have average distance equal to one. Zebrafish \textit{in vivo} tracks and simulated tracks are available in the \texttt{results/species\_tracks} folder in the GitHub repository~(see~\nameref{sec:codeAvail}), together with a notebook to visualise the results and a notebook that shows the pre-processing.

\section{Persistent homology and generator computations}

\subsection{Persistent Homology}

Let $X = \{v_0, \dots, v_M \}$ denote the collection of vertices of a piecewise linear curve. To study the shape of $X$, one can construct a simplicial complex $X_{\varepsilon}$ that represents the connectivity of the points upto a specified distance $\varepsilon$ as the following. Starting with the collection $X$, given a pair of points $v_0, v_1 \in X$ whose distance $d(v_0, v_1) \leq \varepsilon$, we add a 1-simplex between the vertices $v_0$ and $v_1$. Similarly, given a triple of points $v_0, v_1, v_2 \in X$, we add a 2-simplex $[v_0, v_1, v_2]$ if $d(v_0, v_1) \leq \varepsilon$, $d(v_1, v_2) \leq \varepsilon$, and $d(v_0, v_2) \leq \varepsilon$. In general, given an $n$-tuple of points $v_0, \dots, v_{n-1}$, we add a $n$-simplex $[v_0, \dots, v_{n-1}]$ if the distance between every pair of points is at most $\varepsilon$ \footnote{This is often called the Vietoris-Rips complex at parameter $\varepsilon$.}

One can then extract various shape descriptors, such as connected components, loops, and voids, of $X_{\varepsilon}$ by computing homology at various dimensions\footnote{In this paper, we compute homology with field coefficients $\mathbb{F}_2$.}. In particular, the loops of $X_{\varepsilon}$ is encoded by homology in dimension 1, which is indicated by $H_1(X_{\varepsilon})$. Since homology is computed with field coefficients, the dimension-1 homology $H_1(X_{\varepsilon})$ is a vector space, and its dimension indicates the number of loops present in $X_{\varepsilon}$. We refer to specific elements $[x] \in H_1(X_{\varepsilon})$ as \textit{homology classes}. 
 
Recall that $X_{\varepsilon}$ captures the connectivity of points upto distance $\varepsilon$. In order to study the shape of $X$ at various distance resolutions, one can construct the simplicial complexes at various distance parameters $\{ \varepsilon_i \}_{i=0}^N$. This results in the following sequence of simplicial complexes
\begin{equation}
\label{eq:SC_filtration}
X_{\varepsilon_0} \hookrightarrow X_{\varepsilon_1} \hookrightarrow \cdots \hookrightarrow X_{\varepsilon_N}. 
\end{equation}
Note that given two parameters $\varepsilon_i < \varepsilon_j$, there is an inclusion $X_{\varepsilon_i} \hookrightarrow X_{\varepsilon_j}$.

Persistent homology studies the birth and death of various shape descriptors as one varies the parameter $\varepsilon$. To study the evolution of loops, one can compute the homology in dimension 1 for each simplicial complex in Equation \ref{eq:SC_filtration}. This results in the following sequence of vector spaces and linear maps called the persistence module
\begin{equation}
\label{eq:persistence_module}
\mathbb{X}: \quad H_1 (X_{\varepsilon_0}) \to H_1(X_{\varepsilon_1}) \to \dots \to H_1(X_{\varepsilon_N}).
\end{equation} 
In order to extract an interpretable summary of $\mathbb{X}$, one can compute the barcode of $\mathbb{X}$, which is a collection of points $(b_i, d_i)$, where $b_i$ and $d_i$ denote the birth and death parameters of the $i^{\text{th}}$ feature. Structure Theorem \cite{zomorodian2005computing} states that one can always compute such barcode from finite persistence modules of the form in Equation \ref{eq:persistence_module} and that the barcode is an invariant of $\mathbb{X}$. The barcode can be visualised through a persistence diagram, which plots a point for each feature with the $x$-axis representing the birth parameter and the $y$-axis representing the death parameter. 

For each point in the persistence diagram, one can find a \textit{generator}, which is a specific loop that corresponds to the selected point. Finding a generator is a two-step process, both of which requires some choice to be made by the researcher. First, a specific homology class has to be selected, and second, a specific generator has to be selected.  

Let $p$ be a point in the persistence diagram with birth parameter $b$ and death parameter $d$. In the first step, one finds the homology class that corresponds to $p$. That is, one finds the sequence of homology classes $([x_{\epsilon_0}], \cdots [x_{\epsilon_N}])$ such that $[x_{\varepsilon_0}] \in H_1(X_{\varepsilon_0})$, $[x_{\varepsilon_1}] \in H_1(X_{\varepsilon_1})$, $\cdots$, $[x_{\varepsilon_N}] \in H_1(X_{\varepsilon_N})$ and $[x_{\varepsilon}] = 0$ for $\varepsilon \notin [b, d)$. Furthermore, the sequence must be consistent with the maps of Equation \ref{eq:persistence_module}. That is, if $\iota: H_1(X_{\varepsilon_i}) \to H_1(X_{\varepsilon_{i+1}})$ represents the induced map in Equation \ref{eq:persistence_module}, then $[x_{\varepsilon_{i+1}}] = \iota[x_{\varepsilon_i}]$. We refer to this collection $([x_{\epsilon_0}], \cdots [x_{\epsilon_N}])$ as the \emph{persistent homology class}. The persistent homology class of $p$ depends on the choice of basis of the persistence module $\mathbb{X}$.

Once the persistent homology class of $p$ is fixed, one can find the \textit{generator} of $p$. This is usually found by finding the homology class $[x_b] \in H_1(X_b)$ at the birth parameter $b$ and then returning a specific generator of $[x_b]$. For a detailed discussion on homology classes and generators, we refer the reader to \cite{Hatcher}.

\subsection{Minimal generators} 
\label{sec:generators} 
While many persistent homology software return some generator, the generators are often chosen in an aribtrary manner. We present the minimal generators as a principled way of choosing generators. We first compute the generators that minimises the length of the loops. We do this by solving an optimization problem with rational coefficients \cite{LiMinimalCycle} and interpreting the output as generators with $\mathbb{F}_2$ coefficients. We then perform jump-minimization \cite{emmett2015multiscale} to find the generator with the smallest amount of jumps between consecutive points. 

\noindent\textbf{Length minimised generators.} 
Given a persistence diagram, we fix the persistent homology classes using the standard matrix decomposition procedure \cite{VinesVineyardsPersistence}. Note that a generator is a collection of edges that form a loop(s). For each persistent homology class, we find the generator that minimises the sum of the edge lengths. The straightforward approach is to solve the binary optimization problem that minimises the total length of the loop while restricting the coefficients to be ${0,1}$. As this is an NP-hard problem, a variety of combinatorial techniques have been proposed \cite{ChenGenerators, Obayashi2018VolumeOC, Dey2019}. 

In \cite{LiMinimalCycle}, authors approach the NP-hard problem by allowing the solutions to have rational coefficients. This related optimization problem can be solved using linear programming. Their computational experiments show that nearly all solutions resulted in generators with $\{-1, 0, 1\}$ coefficients. Thus, one can interpret the length-minimised generators over the rational coefficients as length-minimised generators with $\mathbb{F}_2$ coefficients. We refer the reader to Algorithm 1 of \cite{LiMinimalCycle} for details.

In our experiments, all length minimization over rational coefficients resulted in generators with coefficients in $\{-1, 0, 1\}$. In theory, it is possible for the length-minimal generators over the $\mathbb{Q}$ to have coefficients in $\{-1, 0, 1\}$. When one encounters such solution, then one could convert all nonzero coefficients to $1$, effectively forcing all coefficients to be in $\mathbb{F}_2$. Note that the outcome isn't necessarily a loop. 
Once we find the length-minimised generator with coefficients in $\mathbb{F}_2$, we denote the collection of vertices that constitute the generator as $(v_0, \dots, v_n)$ with $v_0 < v_1 < \dots < v_n$.

\noindent\textbf{Jump-minimization.}
Given a point $p$ in the persistence diagram with birth parameter $b$ and death parameter $d$, let $(v_0, \dots, v_n)$ denote its length-minimised generator. During this second step, we iteratively update the generator by preferring generators that have smaller jumps between consecutive vertices \cite{emmett2015multiscale}. That is, let $g$ denote the collection of $1$-simplices that represent the generator $(v_0, \dots, v_n)$. Assume that there exists a vertex $v_*$ such that $v_i < v_* < v_{i+1}$, and let $g_*$ denote the colelction of $1$-simplices that represent the generator $(v_0, \dots, v_i, v_*, v_{i+1}, \dots, v_n)$. If $[g] = [g_*]$ in $H_1(X_b)$, then we update the current generator $(v_0, \dots, v_n)$ with $(v_0, \dots, v_i, v_*, v_{i+1}, \dots, v_n)$. Algorithmically, we start with the early vertex $v_0$ and iteratively perform jump-minimization until we reach the last vertex $v_n$. 

Recall that there may be rare cases in which the length-minimization over the rational coefficients isn't necessarily a loop in $H_1(X_b)$ with $\mathbb{F}_2$ coefficients. In such cases, we still perform jump-minimization by testing whether the 2-simplex $[v_i, v_*, v_{i+1}]$ is trivial in $H_1(X_b)$ or not. 

\begin{figure}[!ht]
    \includegraphics[width = 10cm]{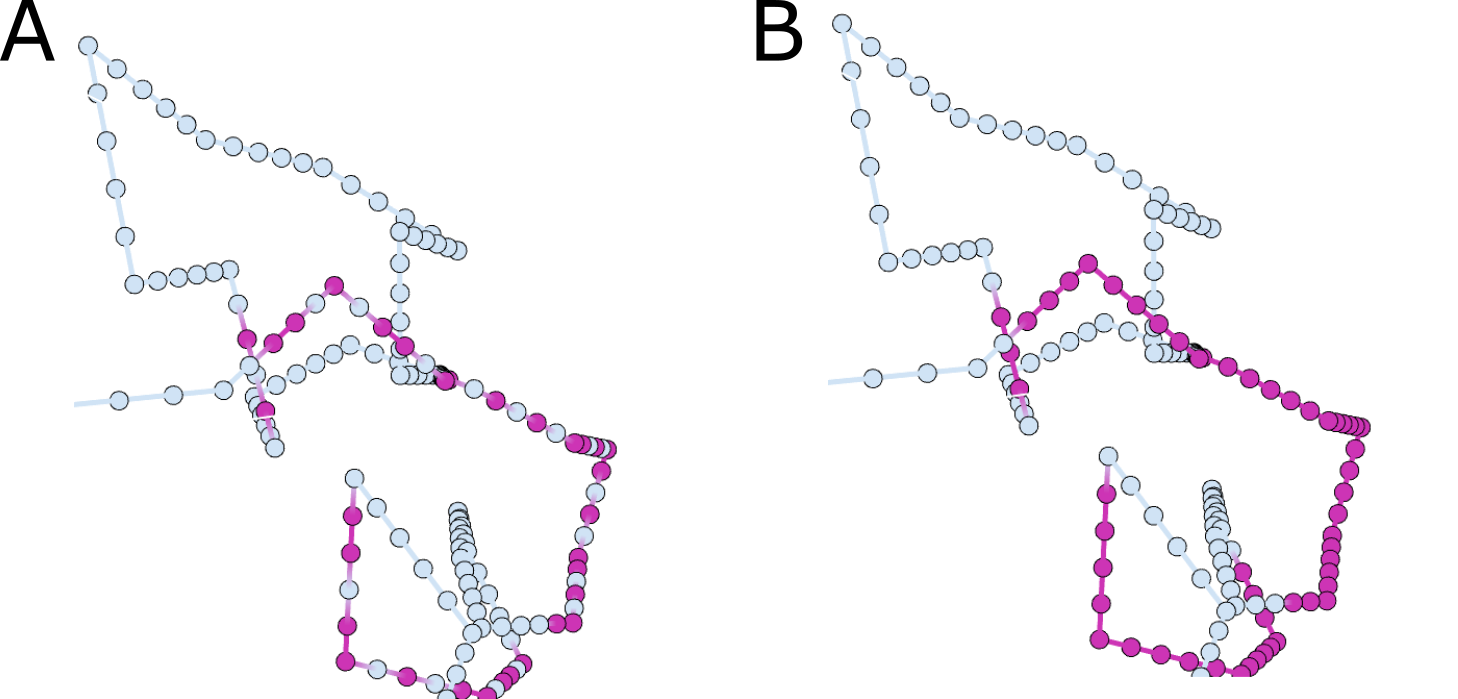}
    \captionsetup{margin=-0.8cm}

    \caption{\textbf{A.} Example length minimised generator. \textbf{B.} Generator after jump minimization. Given a collection of points aligned in a line segment, the length-minimised generator consists of an arbitrary subset of points. The jump minimization ensures that all points in a given line segment are included in the generator.}
    \label{fig:generator_steps}
\end{figure}

Figure~\ref{fig:generator_steps} illustrates the need for such jump minimization step. Figure~\ref{fig:generator_steps}A illustrates an example length-minimal generator over the rationals. When one has a collection of points that are in a line segment, the minimal generator will make arbitrary choices in terms of which points are included in the minimal generator. The jump minimization step  Figure~\ref{fig:generator_steps}B.

The code for computing minimal generators is available in the minimal generators GitHub repository~(see~\nameref{sec:codeAvail}).

\section{Analysis}

\subsection{PH-community analysis}
To characterise spatial features and find structural differences in trajectories, we analyse the communities induced by hyperTDA. Given a curve C and its corresponding hyperTDA-induced partition into communities, for each community $c$ we compute:

\begin{itemize}
    \item The \textit{volume} of the sub-curve spanned by vertices in the PH-community is quantified as its \textit{radius of gyration} (ROG) and computed using the \texttt{spacecurves.radius\_of\_gyration()} function of the \texttt{pyknotid} \cite{pyknotid} Python software.
    
    \item The \textit{PH-community size}. This is computed as the L$^1$ norm of the corresponding (possibly \textit{un}interpolated) community vector, \textit{i.e.}~the vector of length equal to the number of points in the curve, whose $i^{th}$ entry is 0 if the point \textit{i} is not a member of \textit{c}, and is equal to 1 (or to a value $v$ with $0 < v \leq 1$ for \textit{un}interpolated communities) otherwise. 
    
    \item The \textit{geodesic size} of the community. This is computed as the difference between the smallest and largest non-zero indices in the community vector. 
    
\end{itemize}

Further, for each pair of communities in a partition we compute their \textit{geodesic intersection}, which is given by the $p$-value of the Mann-Whitney U-test of the corresponding community vectors. Details on this analysis can be found in the notebook \texttt{examples/with\_interpolation/result\_visualisation.ipynb} in the GitHub repository~(see~\nameref{sec:codeAvail}).

\subsection{CNN architecture}
A CNN architecture was designed to distinguish each of the diffusion models, given observations of trajectories with unknown $\alpha$.
The CNN architecture was chosen by initially implementing simple proof-of-concept architectures and adding the least amount of complexity capable of capturing hypothesised patterns in the data. The code for all tested models is available at \texttt{src/hypergraph\_CNN.jl} in the GitHub repository~(see~\nameref{sec:codeAvail}).
The chosen architecture and its performance are shown in Figure~\ref{fig:CNN}. The CNN takes as input the data from a weighted hypergraph $\mathcal{H} = (\mathcal{N},\mathcal{E}, \mathcal{W} )$, equipped with node-weights $\mathcal{V}$ for each trajectory. Specifically, we feed a matrix $H$ and a vector $\mathcal{V}$ (and potentially $\mathcal{W}$) from each trajectory to the architecture. The matrix $H$ has a column for each vertex in $\mathcal{N}$, with entries indicating membership in each PH generator (\textit{i.e.}~each hyperedge in $\mathcal{E}$), making $H$ a hypergraph incidence matrix. In the reported computations, $\mathcal{V}$ is given by node-centrality, while $\mathcal{W}$ (given by persistence) is not used as input, as it does not improve the results. The network was designed to capture patterns at any location along the trajectory, hence convolutions run along the ordered vertices in the first layer, both on the (possibly weighted) $H$ and on $\mathcal{V}$ to separately capture first order patterns in both. Naturally, the kernels only span a single hyperedge since hyperedges are independent vertex collections, with arbitrary ordering. Also for this reason, the resulting tensors (containing first order patterns) are reduced by maxpooling. All hidden values are then concatenated, before a second and final convolution along vertices, which merely serves to combine patterns from V and H, hence the kernel only spans a single vertex. At this point there is no longer a use in distinguishing locations along the trajectory, so the hidden vertices are maxpooled, before being reshaped to the output prediction vector by a simple dense layer.

\section{Additional results}

\subsection{Results with alternative generators}

In this section we show the result of hyperTDA analysis with the alternative (minimal or matroid, depending on what was presented on the main manuscript) on the datasets considered.

\begin{figure}[!ht]
\centering

    \includegraphics[width = 16cm]{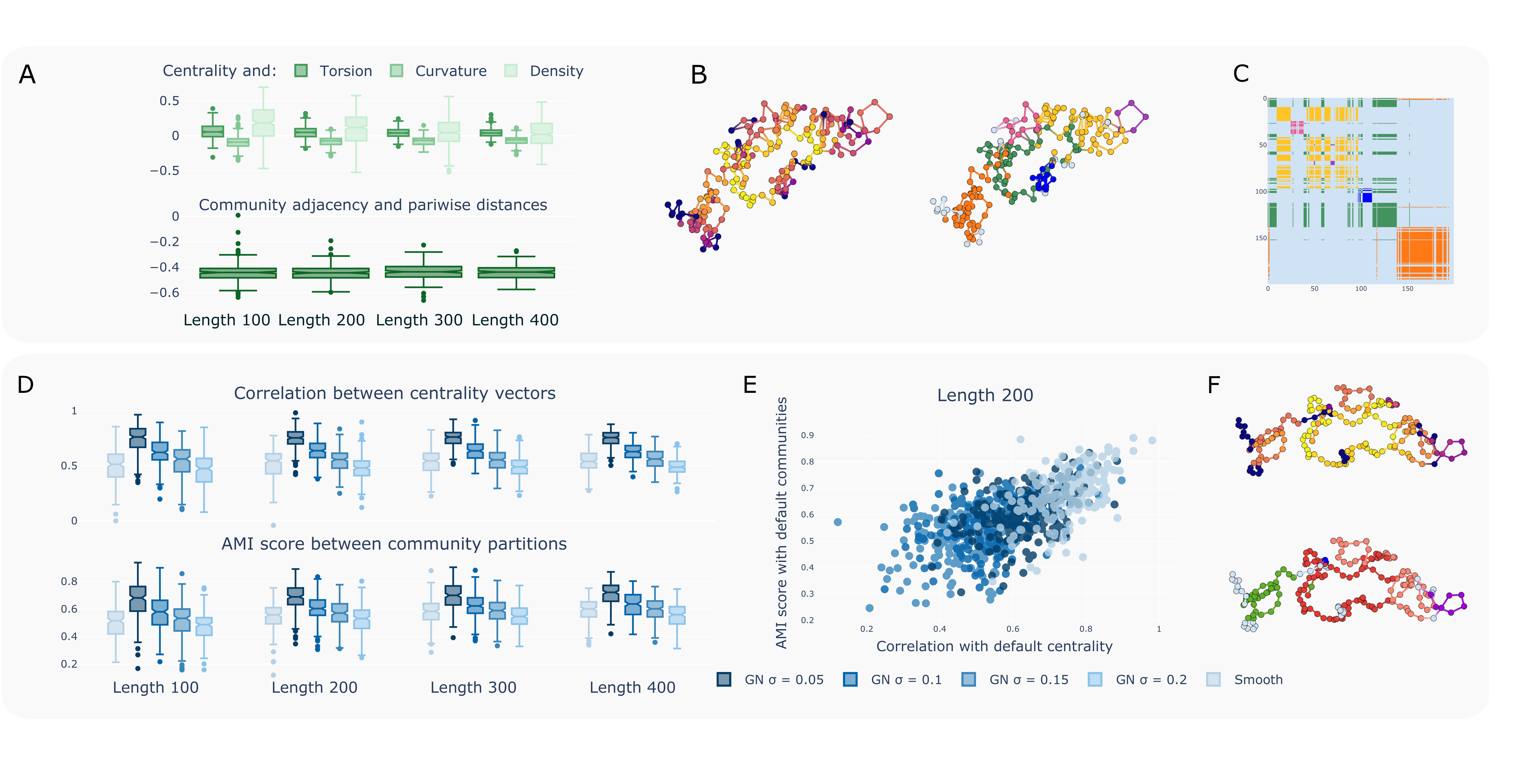}
    \captionsetup{margin=-0.8cm}
    \caption{\textbf{HyperTDA interpretation and robustness with minimal generators.} \textbf{(A)} (Top) Correlation between vertex centrality vectors and curvature, torsion and density vectors. (Bottom) Correlation between PH-community matrices and pairwise distance matrices. \textbf{(B)} Example curve C of length 200 plotted with vertices coloured according to centrality values (top) and to community membership (bottom).  \textbf{(C)} Heatmap of C's PH-community matrix, where non-zero entries are coloured according to community membership. \textbf{(D)} (Top) Correlation between centrality of default curves and noisy ones. (Bottom) Adjusted mutual information (AMI) score between community partitions of default curves and noisy ones. \textbf{(E)} AMI score with default communities plotted against correlation with default centrality, for the length 200 dataset. \textbf{(F)} The curve obtained by smoothing C, plotted with vertices coloured according  to centrality values(top) and to community membership  (bottom). HyperTDA of the smooth curve is very similar to that of the original curve (panel C), indicating the robustness of hyperTDA. In particular, the central loop is still identified as the predominant feature. Out of the 5000 trajectories analysed, 53 had trivial PH and further 15 had minimal generators with non-integers coefficients. These trajectories were removed from the analysis.  }
    \label{fig:validation1_minimal}
\end{figure}

\begin{figure}[!ht]
    \hspace{-7mm}
    \includegraphics[width = 14cm]{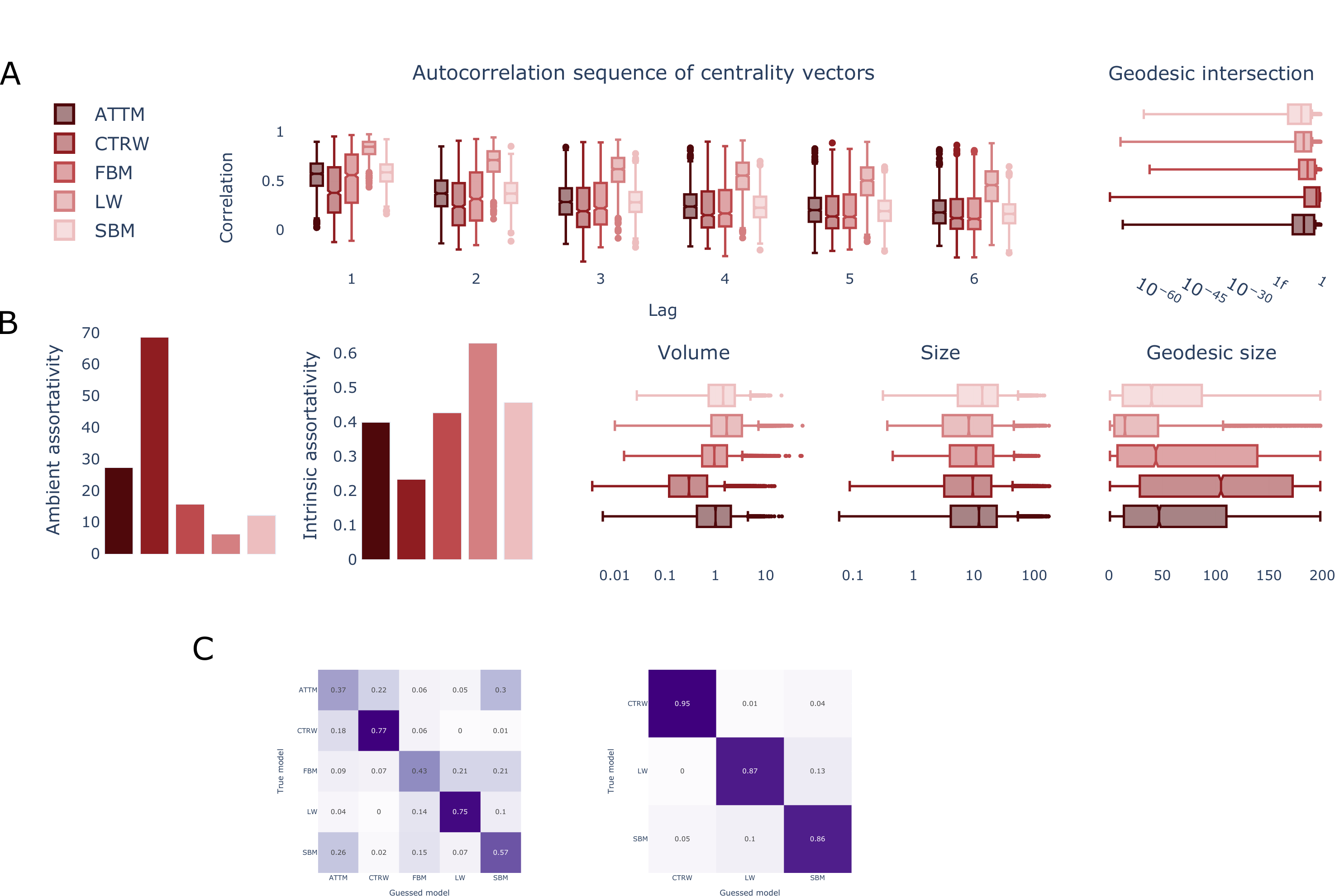}
    \captionsetup{margin=-0.8cm}
    \caption{\textbf{Analysis of AnDi trajectories with minimal generators}
    \textbf{(A)} Auto-correlation of centrality distributions, by model; for each value $m$ in $\{1,\cdots{}, 6\}$ the boxplots show the distribution of correlation values between centrality vectors and their values shifted by a lag $m$. \textbf{(B)} Analysis of the trajectories' communities, by model. The left-most plot shows the bar plot given by model ambient and intrinsic assortativities. Box plots for the distributions of volume, size, geodesic size and geodesic intersection are shown on the right. \textbf{(C)} The confusion matrix summarising the classification accuracy between all the 5 models, and the one summarising the classification accuracy between CTRW, LW and SBM models. }
    \label{fig:SIAnDi}
\end{figure}

\begin{figure}[!ht]
    \includegraphics[width = 10cm]{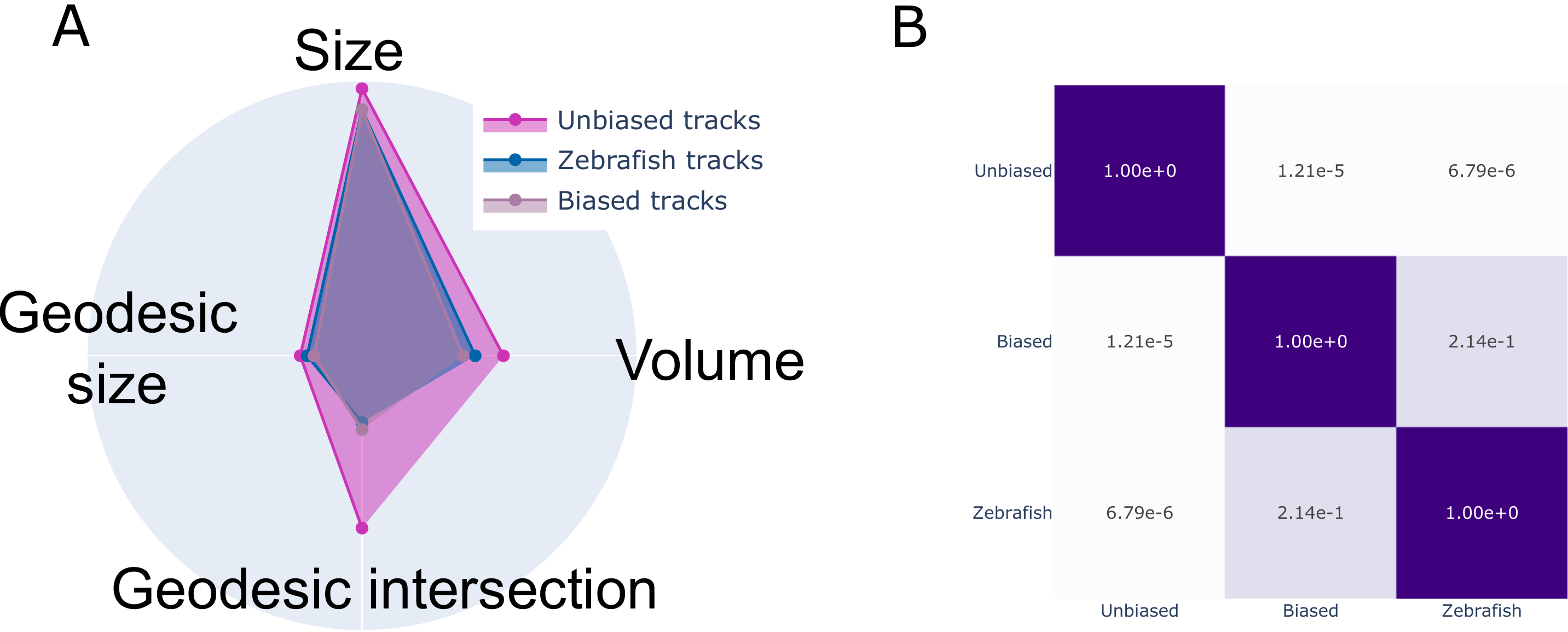}
    \captionsetup{margin=-0.8cm}

    \caption{\textbf{Cell migration trajectories with matroid generators} \textbf{(A)}  Radar chart given by the medians of volume, size, geodesic size and pairwise geodesic intersection distributions for the unbiased and biased simulated data, and for the zebrafish trajectories. \textbf{(B)} Heat-map of pairwise Kolmogorov-Smirnov test on the cumulative distributions obtained by concatenating geodesic intersection, volume, size and geodesic size distributions.}
    \label{fig:animalcellSI}
\end{figure}

\begin{figure}[!ht]
    \includegraphics[width = 9cm]{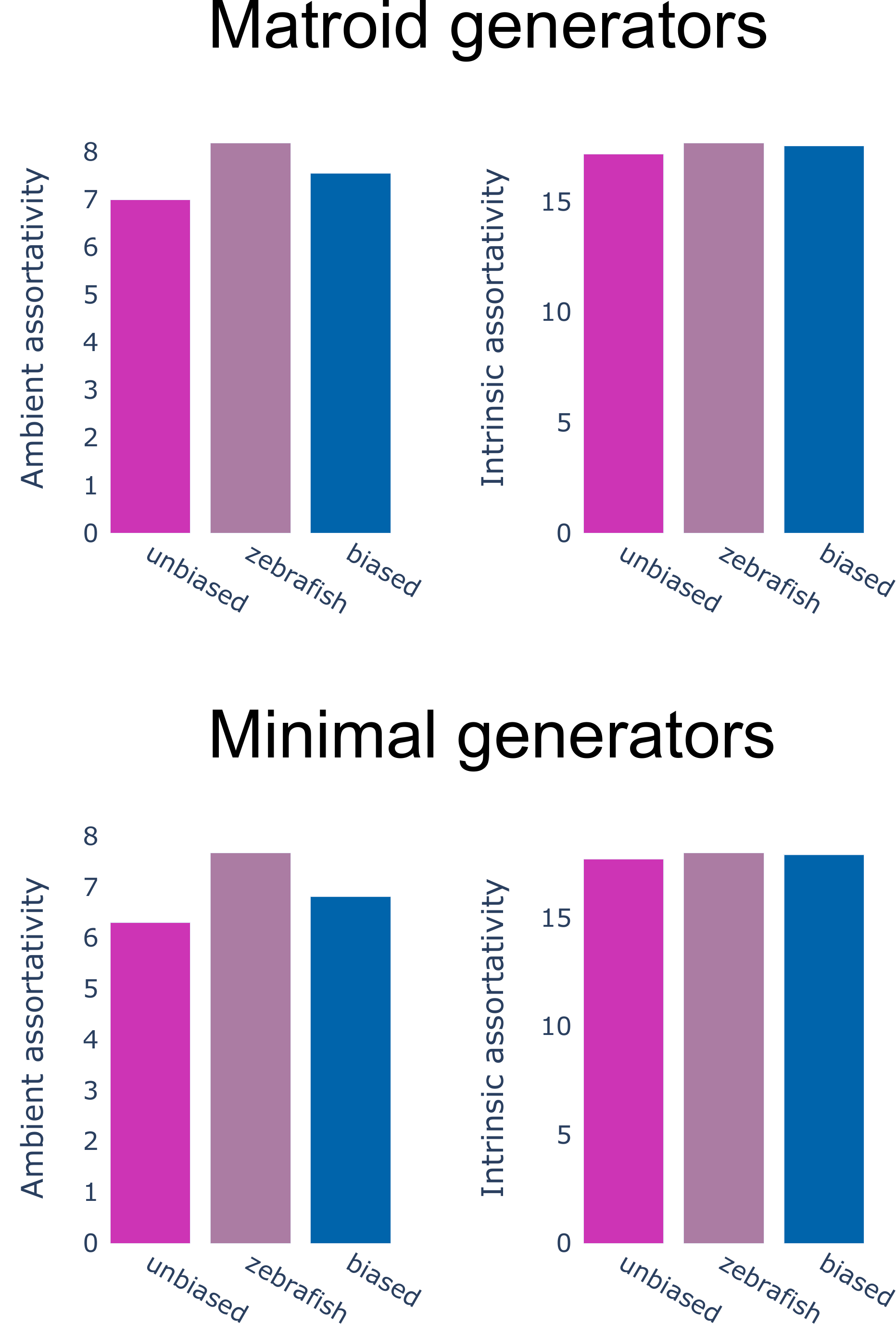}
    \captionsetup{margin=-0.8cm}

    \caption{\textbf{Assortativities for zebrafish, biased and unbiased trajectories with matroid and minimal generators}}
    \label{fig:animalcellSI2}
\end{figure}

\begin{figure}[!ht]
    \includegraphics[width = 14cm]{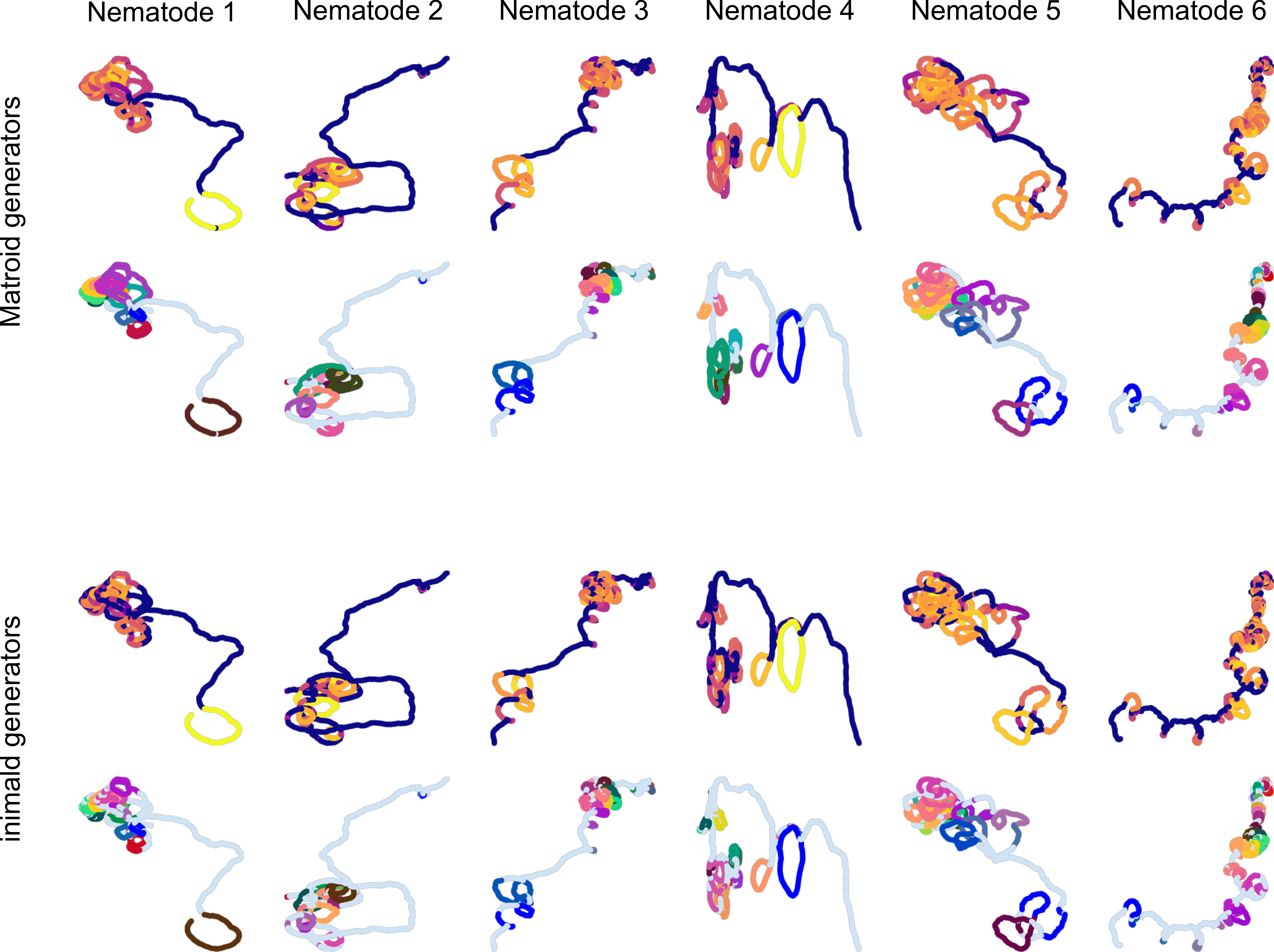}
    \captionsetup{margin=-0.8cm}

    \caption{\textbf{Nematodes tracks with matroid and minimal generators} From top to bottom: Nematode tracks with vertices coloured by centrality values and computed with matroid generators.  Nematode tracks with vertices coloured by community membership and computed with matroid generators.  Nematode tracks with vertices coloured by centrality values and computed with minimal generators. Nematode tracks with vertices coloured by community membership and computed with minimal generators.}
    \label{fig:animalcellSI3}
\end{figure}

\begin{figure}[!ht]
\centering

    \includegraphics[width = 16cm]{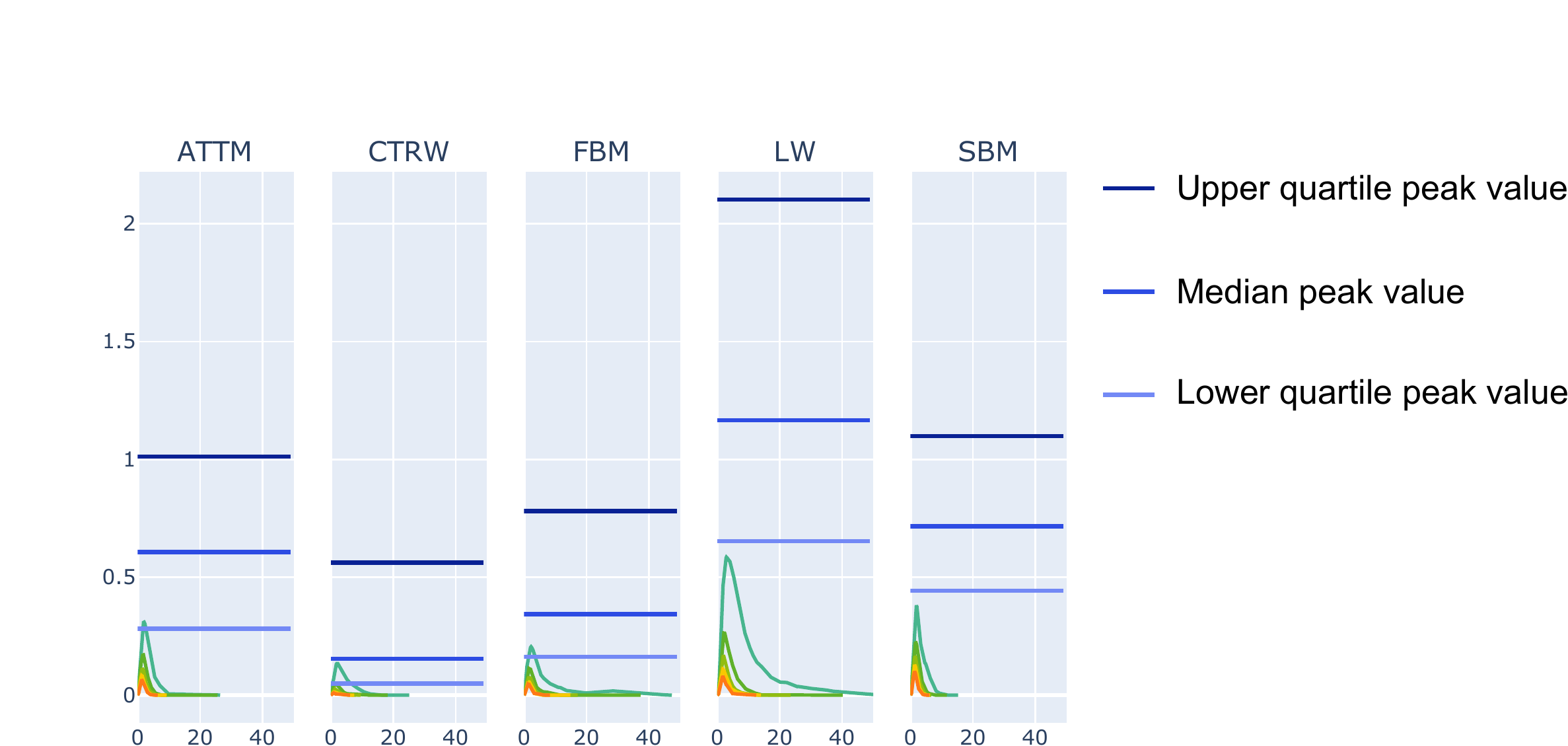}
    \captionsetup{margin=-0.8cm}
    \caption{\textbf{Average persistent landscapes of AnDi trajectories.} Persistent landscapes are an alternative representation of PH's output, equivalent to persistent diagrams. Landscapes in the figure are computed with the code implemented in \cite{benjamin2022homology}. The figure shows the persistent landscapes obtained by averaging the persistent landscapes of AnDi trajectories, by model. Horizontal lines in the plots show the median, lower and upper quartiles of the $\lambda_1$'s max values, by model. The average landscapes look different. However, unlike in the PH-community analysis, those differences are not directly interpretable as model-specific structural features. Further, even though the averages are not trivial, the information retained is strongly reduced with respect to individual landscapes/diagrams. This is shown by the lack of structure (the average landscapes are just nested peaks localised towards the origin of the x-axis). Further, the peak values of each $\lambda$ function are much lower than in individual landscapes, as shown by the median, lower and upper quartiles of  $\lambda_1$'s peaks. }
    \label{fig:landscapes}
\end{figure}

\end{document}